\documentclass[a4paper,12pt]{article}

\usepackage{amsmath,amssymb,amsthm,amscd,cases}

{}
\numberwithin{equation}{section}

\newtheorem{thm}{Theorem}[section]

\newcommand{\R}{\mathbb{R}}

\renewcommand{\a}{\alpha}
\renewcommand{\b}{\beta}
\renewcommand{\d}{\delta}
\newcommand{\p}{\partial}

\newcommand{\bc}{\begin{cor}}
\newcommand{\ec}{\end{cor}}
\newcommand{\bl}{\begin{lem}}
\newcommand{\el}{\end{lem}}
\newcommand{\bp}{\begin{prop}}
\newcommand{\ep}{\end{prop}}
\newcommand{\bt}{\begin{thm}}
\newcommand{\et}{\end{thm}}
\newcommand{\bal}{\begin{array}{ll}}
\newcommand{\ba}{\begin{array}}
\newcommand{\bac}{\begin{array}{ccc}}
\newcommand{\ea}{\end{array}}

\newcommand{\be}{\begin{equation}}
\newcommand{\ee}{\end{equation}}

\newtheorem{lem}[thm]{Lemma}
\newtheorem{cor}[thm]{Corollary}
\newtheorem{prop}[thm]{Proposition}

\begin{document}

\title{\Large\bf Global Existence and Asymptotic Behavior \\
of Solutions for Quasi-linear Dissipative \\
Plate Equation}

\author{Yongqin Liu{\footnote{email:~yqliu2@yahoo.com.cn}}\\
 \small\emph{Faculty of Mathematics, Kyushu University}\\
 \small\emph{Fukuoka 819-0395, Japan} \\[3mm]
Shuichi Kawashima{\footnote{email:~kawashim@math.kyushu-u.ac.jp}}\\
 \small\emph{Faculty of Mathematics, Kyushu University}\\
 \small\emph{Fukuoka 819-0395, Japan}
 }

\date{}

\maketitle



\begin{abstract}

In this paper we focus on the initial value problem for quasi-linear 
dissipative plate equation in multi-dimensional space $(n\geq2)$. 
This equation verifies the decay property of the regularity-loss type, 
which causes the difficulty in deriving the global a priori estimates 
of solutions. We overcome this difficulty by employing a 
time-weighted $L^2$ energy method which makes use of 
the integrability of $\|(\p^2_xu_t,\p^3_xu)(t)\|_{L^{\infty}}$. 
This $L^\infty$ norm can be controlled by showing the optimal $L^2$ 
decay estimates for lower-order derivatives of solutions. 
Thus we obtain the desired a priori estimate which enables us to 
prove the global existence and asymptotic decay of solutions 
under smallness and enough regularity assumptions on the initial data. 
Moreover, we show that the solution can be approximated by 
a simple-looking function, which is given explicitly in terms of 
the fundamental solution of a fourth-order linear parabolic equation. 
\\



\end{abstract}


\section{Introduction}

In this paper we consider the initial value problem of the
following quasi-linear dissipative plate equation in 
multi-dimensional space ${\mathbb R}^n$ with $n\geq 2$: 
\begin{equation}\label{1a}
u_{tt}-\Delta u_{tt}
+\sum_{i,j=1}^n b^{ij}(\partial_x^2u)_{x_ix_j}+u_t=0.
\end{equation}
The initial data are given as 
\begin{equation}\label{IC}
u(x,0)=u_0(x), \quad u_t(x,0)=u_1(x).
\end{equation}
Here $u=u(x,t)$ is the unknown function of 
$x=(x_1,\cdots, x_n)\in\R^n$ and $t>0$, and represents the 
transversal displacement of the plate at the point $x$ and 
the time $t$. The term $u_t$ represents a frictional dissipation 
in the plate, and the term $\Delta u_{tt}$ corresponds to 
the rotational inertia effects. 
Also, $b^{ij}=b^{ij}(V)$ are smooth functions of 
$V=(V_{ij})\in{\mathcal S}_n$ satisfying the following structural 
conditions, where ${\mathcal S}_n$ denotes the totality of 
$n\times n$ real symmetric matrices: 

\begin{enumerate}

\item
[{[A1]}]\ 
There is a smooth potential function $\phi=\phi(V)$ such that 
$b^{ij}(V)=(\partial\phi/\partial V_{ij})(V)$ for 
$V=(V_{ij})\in{\mathcal S}_n$. 
\item
[{[A2]}]\ 
$\gamma(\omega):=\sum_{ij}\sum_{\alpha\beta}b^{ij}_{\alpha\beta}(O)
\omega_i\omega_j\omega_\alpha\omega_\beta>0$ for 
$\omega=(\omega_1,\cdots,\omega_n)\in S^{n-1}$, 
where the summations are taken over all $i,j=1,\cdots,n$ and 
$\alpha,\beta=1,\cdots,n$, respectively, and $O$ denotes the zero 
matrix in ${\mathcal S}_n$. 
\end{enumerate}

\noindent
Here we put 
\begin{equation}\label{1.3}
b^{ij}_{\alpha\beta}(V)
=\frac{\partial b^{ij}}{\partial V_{\alpha\beta}}(V)
=\frac{\partial^2\phi}
{\partial V_{ij}\partial V_{\alpha\beta}}(V), \quad
i,j,\alpha,\beta=1,\cdots,n.
\end{equation}
By the definition we see that 
\begin{equation}\label{1.4}
b^{ji}_{\alpha\beta}(V)=b^{ij}_{\alpha\beta}(V)
=b^{ij}_{\beta\alpha}(V), \qquad
b^{ij}_{\alpha\beta}(V)=b^{\alpha\beta}_{ij}(V)
\end{equation}
for $V\in{\mathcal S}_n$ and $i,j,\alpha,\beta=1,\cdots,n$. 
Namely, $b^{ji}_{\alpha\beta}(V)$ is symmetric with respect 
to the indexes $(i,j)$ and $(\alpha,\beta)$. 
We note that the condition [A2] implies 
\begin{equation}\label{elliptic}
\sum_{i,j=1}^n\sum_{\alpha,\beta=1}^n\int_{{\mathbb R}^n}
b^{ij}_{\alpha\beta}(O)u_{x_ix_j}u_{x_\alpha x_\beta}dx
\geq c\|\partial_x^2 u\|_{L^2}^2
\end{equation}
with a positive constant $c$. 
The potential function $\phi(V)$ in [A1] is called the free energy 
of \eqref{1a}. Throughout the paper, without loss of generality, 
we may assume that 
\begin{equation}\label{1.6}
\phi(O)=0, \qquad b^{ij}(O)=0, \quad i,j=1,\cdots,n.
\end{equation}
In this case we have 
$b^{ij}(V)=\sum_{\alpha\beta}b^{ij}_{\alpha\beta}(O)V_{\alpha\beta}
+O(|V|^2)$, so that the linearized equation of \eqref{1a} becomes to 
\be\label{1cc}
u_{tt}-\Delta u_{tt}+\sum_{i,j=1}^n\sum_{\alpha,\beta=1}^n
b^{ij}_{\alpha\beta}(O)u_{x_ix_jx_\alpha x_\beta}+u_t=0.
\ee

Equations of the fourth-order appear in problems of solid mechanics 
and in the theory of thin plates and beams, 
and elliptic equations of the fourth-order appear in some formulations 
of problems related to the Navier-Stokes equations (see \cite {Te}). 
L. S\'{a}nchez \cite{Sa} studied the existence and uniqueness of 
solutions of the full von K\'{a}rm\'{a}n system in an exterior domain 
and in the whole space, and proved that the model for thermoelastic 
plates is a singular limit of the von K\'{a}rm\'{a}n system under 
thermal effects. 
In \cite{MZ} Perla Menzala and Zuazua also showed that the plate 
equation can be obtained as a singuler limit of the von K\'{a}rm\'{a}n 
system. 
Enomoto \cite{E} studied a linear thermoelastic system associated 
with the plate equation in an exterior domain and proved a polynomial 
decay of the local energy for initial data with compact support. 
Buriol \cite{Bu} considered the Timoshenko system of thermoelastic 
plates in $\R^n$ and showed an exponential decay of the energy by 
considering two types of dissipation in the system. 
In \cite{CBBP}, Char$\tilde{\rm a}$o {\it et al.} studied 
the asymptotic behavior of solutions of a dissipative plate equation 
in $\R^n$ with periodic coefficients. They used the Bloch waves 
decomposition and a convenient Lyapunov function to derive a complete 
asymptotic expansion of solutions as $t\to\infty$, and also proved 
that the solutions for the linear model behave as the homogenized 
heat kernel. 

In \cite{CL}, da-Luz and Char$\tilde{\rm a}$o studied a semi-linear 
dissipative plate equation whose linear part is given by \eqref{1cc} with 
$b^{ij}_{\alpha\beta}(O)=\delta_{ij}\delta_{\alpha\beta}$:
\begin{equation}\label{1c}
u_{tt}-\Delta u_{tt}+\Delta^2u+u_t=0.
\end{equation}
They proved the global existence of solutions and a polynomial 
decay of the energy by exploiting an energy method. 
However their result was restricted to the lower dimensional case 
$1\leq n \leq5$. 
This restriction on the space dimension was removed by 
Sugitani and Kawashima \cite{SK} by making use of the sharp decay 
estimates for the linearized equation \eqref{1c}. 
One of the decay estimates obtained in \cite{SK} for \eqref{1c} is 
given as follows: 
When $u_0\in H^{s+1}\cap L^1$ and $u_1\in H^s\cap L^1$, 
\begin{equation*}\label{}
\|\p^k_xu(t)\|_{H^{s-\sigma_1(k,n)}}\leq
C(1+t)^{-{n\over8}-{k\over4}},
\end{equation*}
where $s\geq\sigma_1(k,n):=k+[{{n+2k-1}\over4}]$. 
This implies that we have the optimal decay rate 
$t^{-{n\over8}-{k\over4}}$ for $\p^k_xu$ only by assuming 
the additional $\sigma$-th order regularity on the initial data, 
where $\sigma:=\sigma_1(k,n)-k=[{{n+2k-1}\over4}]$. 
This is the decay structure of the regularity-loss type and is 
characterized by the property 
$$
{\rm Re}\lambda(\xi)\leq -c|\xi|^4/(1+|\xi|^2)^3, 
$$
where $\lambda(\xi)$ denotes the eigenvalue of the equation obtained 
by taking Fourier transform of \eqref{1c}. 
See \cite{SK} for the details. 
A similar decay structure of the regularity-loss type was also 
observed for the dissipative Timoshenko system (\cite{IK,RR}) and 
a hyperbolic-elliptic system related to a radiating gas (\cite{HK}). 
For more studies on various aspects of dissipation of 
plate equations, we refer to \cite{BL,DK,Le,LZ2,MNV,Pa,PVM,ZZ}. 
Also, as for the study of decay properties for dissipative 
hyperbolic-type equations, we refer to \cite{LW,Ma,Ni,TY} for 
damped wave equations and \cite{FL,LZ1,Mu} for wave equations 
of memory-type dissipation. 

The main purpose of this paper is to 
study the global existence and large-time behavior of 
solutions to the initial value problem \eqref{1a}, \eqref{IC}. 
For this problem, we will observe that the regularity-loss also 
occurs in the dissipative part of the energy estimates, 
which makes it difficult to show the global a priori estimates 
of solutions to the nonlinear problem. 
To overcome the difficulty, we employ a time-weighted 
energy method, in which we make good use of the integrability of 
$\|(\p^2_xu_t,\p^3_xu)(t)\|_{L^{\infty}}$ and the optimal decay of 
the lower-order derivatives of solutions. 
Consequently, we obtain the global existence and the optimal decay 
estimates of solutions for small initial data with sufficient 
regularity. Moreover, we also show that for $t\to\infty$, 
the global solution obtained is asymptotic to a simple-looking 
function which is given explicitly in terms of the the fundamental 
solution of the fourth-order linear parabolic equation 
\begin{equation}\label{limit-eq}
u_t+\sum_{i,j=1}^n\sum_{\alpha,\beta=1}^n
b^{ij}_{\alpha\beta}(O)u_{x_ix_jx_\alpha x_\beta}=0.
\end{equation}

The contents of the paper are as follows. In Section 2 we give full 
statements of our main theorems. A time-weighted energy method is 
introduced in Section 3, and we explain why the usual energy method 
does not work well for our problem. In Section 4, we give the 
optimal decay estimates of solutions. Section 5 gives the proof of 
the first main theorem on the global existence and the optimal decay 
estimates of solutions, and the last section gives the proof of 
the second theorem on the asymptotic profile. 

Before closing this section, we give some notations to be used 
below. Let $\mathcal{F}[f]$ denote the Fourier transform of $f$ 
defined by 
$$
\mathcal{F}[f](\xi)=\hat{f}(\xi)
:=\int_{\R^n}e^{-ix\cdot\xi}f(x)dx,
$$
and we denote its inverse transform by $\mathcal{F}^{-1}$. 

For $1\leq p\leq\infty$, $L^p=L^p(\R^n)$ is the usual Lebesgue space 
with the norm $\|\cdot\|_{L^p}$. 
For $\gamma\in\R$, let $L^1_{\gamma}=L^1_{\gamma}(\R^n)$ denote the
weighted $L^1$ space with the norm 
$$
\|f\|_{L^1_{\gamma}}:=\int_{\R^n}(1+|x|)^{\gamma}|f(x)|dx.
$$
Let $s$ be a nonnegative integer. 
Then $H^s=H^s(\R^n)$ denotes the Sobolev space of $L^2$ functions, 
equipped with the norm 
$$
\|f\|_{H^s}
:=\Big(\sum\limits_{k=0}\limits^{s}
\|\partial_x^kf\|_{L^2}^2\Big)^{1\over2}.
$$
Here, for a nonnegative integer $k$, $\partial_x^k$ denotes the totality 
of all the $k$-th order derivatives with respect to 
$x\in{\mathbb R}^n$. 
Also, $C^k(I; H^s(\R^n))$ denotes the space of $k$-times 
continuously differentiable functions on the interval $I$ 
with values in the Sobolev space $H^s=H^s(\R^n)$. 

Finally, in this paper, we denote every positive constant by the 
same symbol $C$ or $c$ without confusion. 
$[\,\cdot\,]$ is Gauss' symbol.


\section{Main theorems}

Our first theorem is on the global existence and optimal decay of 
solutions to the problem \eqref{1a}, \eqref{IC}.
To state the result, we need to introduce several special notations. 
Let $k$ be a nonnegative integer and $n$ be the space dimension. 
Let 
\begin{eqnarray}\label{2.1}
&\sigma_0(k)=k+[{{k+1}\over2}], \qquad
\sigma_1(k,n)=k+[{{n+2k-1}\over4}], \\[-6mm]
& \nonumber
\end{eqnarray}
and put $\sigma(k,n)={\rm max}\{\sigma_0(k),\,\sigma_1(k,n)\}$, 
which indicates the loss of regularity. 
Since $\sigma_1(k,3)=\sigma_0(k)$ and $\sigma_1(k,n)$ is an 
increasing function of $n$, we have 
\begin{equation}\label{2.2}
\sigma(k,n)=\left\{\bal 
\sigma_0(k), \ & 1\leq n\leq 3,\\[1mm]
\sigma_1(k,n), \ & n\geq 3.
\ea\right.
\end{equation}
Also, for $n\geq 2$, we define $s(n)$ by 
\begin{equation}
s(n)=\left\{\bac 
8, \ & n=2,\\
6, \ & n=3,\\
3[{n\over4}]+5,\ &n\geq 4,
\ea\right.
\end{equation}
which indicates the regularity of the initial data. 
Now we can state our first theorem as follows. 

\bt[Global existence and optimal decay]\label{21} 
Suppose that the conditions {\rm [A1]} and {\rm [A2]} are 
satisfied. Let $n\geq 2$ and $s\geq s(n)$. 
Assume that $u_0\in H^{s+1}(\R^n)\cap L^1(\R^n)$ and 
$u_1\in H^s(\R^n)\cap L^1(\R^n)$, and put 
$$
E_1:=\|u_0\|_{H^{s+1}}+\|u_1\|_{H^s}+\|(u_0,u_1)\|_{L^1}.
$$
Then there is a positive constant $\d_0$ such that if 
$E_1\leq\d_0$, then the problem \eqref{1a}, \eqref{IC} has a unique 
global solution $u(x,t)$ with 
$$
u\in C^0([0,\infty); H^{s+1}(\R^n))
\cap C^1([0,\infty); H^s(\R^n)).
$$
Moreover, the solution satisfies the following optimal decay 
estimates: 
\begin{equation}\label{decay-1}
\|\p^k_xu(t)\|_{H^{s-1-\sigma(k,n)}}\leq
CE_1(1+t)^{-{n\over8}-{k\over4}},
\end{equation}
\begin{equation}\label{decay-2}
\|\p^k_xu_t(t)\|_{H^{s-4-\sigma(k,n)}}\leq
CE_1(1+t)^{-{n\over8}-{k\over4}-1}
\end{equation}
for $k\geq 0$, where $\sigma(k,n)\leq s-1$ in \eqref{decay-1} and 
$\sigma(k,n)\leq s-4$ in \eqref{decay-2}. 
\et

\noindent
{\bf Remark.}\ 
The regularity assumption $s\geq s(n)$ might be technical but it seems 
necessary in our proof. 

\medskip
This global existence and optimal decay result is based on the 
following local existence result and the corresponding a priori 
estimates stated in Proposition \ref{AP} below. 

\bt[Local existence]\label{LS}
Suppose that the conditions {\rm [A1]} and {\rm [A2]} are 
satisfied. 
Let $n\geq 1$ and $s\geq [{n\over2}]+3$. 
Assume that $u_0\in H^{s+1}(\R^n)$ and $u_1\in H^s(\R^n)$, and put 
$E_0:=\|u_0\|_{H^{s+1}}+\|u_1\|_{H^s}$. 
Then there is a positive constant $T_0$ depending on $E_0$ 
such that the problem \eqref{1a}, \eqref{IC} has a unique solution 
$u(x,t)$ with 
$$
u\in C^0([0,T_0]; H^{s+1}(\R^n))\cap C^1([0,T_0]; H^s(\R^n)). 
$$
The solution verifies the following estimate for $t\in[0,T_0]$: 
$$
\|u(t)\|_{H^{s+1}}^2+\|u_t(t)\|_{H^s}^2\leq CE_0^2.
$$
\et

This local existence result can be proved by the standard method 
based on the successive approximation sequence, so that the details 
are omitted. 

To state the result on our a priori estimates, we introduce several 
time-weighted norms: 
\begin{equation}\label{norm-1}
\begin{split}
E(T)^2
&:=\sum\limits^{[{{s+1}\over3}]}\limits_{j=0}
  \sup\limits_{0\leq t\leq T}(1+t)^j
  \|\p^{2j}_xu(t)\|^2_{H^{s-3j+1}}
  +\sup\limits_{0\leq t\leq T}\|u_t(t)\|^2_{H^s} \\
&\qquad
  +\sum\limits^{[{{s-2}\over3}]}\limits_{j=0}
  \sup\limits_{0\leq t\leq T}(1+t)^{j+1}
  \|\p^{2j}_xu_t(t)\|^2_{H^{s-3j-1}}, \\
D(T)^2
&:=\sum\limits^{[{{s+1}\over3}]}\limits_{j=1}
  \int^T_0(1+\tau)^{j-1}\|\p^{2j}_xu(\tau)\|^2_{H^{s-3j+1}}d\tau
  +\int^T_0\|u_t(\tau)\|^2_{H^{s-1}}d\tau \\
&\qquad
  +\sum\limits^{[{{s-2}\over3}]}\limits_{j=0}
  \int^T_0(1+\tau)^{j+1}\|\p^{2j}_xu_t(\tau)\|^2_{H^{s-3j-2}}d\tau,
\end{split}
\end{equation}
where $s\geq 2$, and 
\begin{equation}\label{norm-2}
\begin{split}
&M_0(T):=\sum\limits_{\sigma(k,n)\leq s-1}
  \sup\limits_{0\leq t\leq T}(1+t)^{{n\over8}+{k\over4}}
  \|\p^k_xu(t)\|_{H^{s-1-\sigma(k,n)}}, \\
&M_1(T):=\sum\limits_{\sigma(k,n)\leq s-4}
  \sup\limits_{0\leq t\leq T}(1+t)^{{n\over8}+{k\over4}+1}
  \|\p^k_xu_t(t)\|_{H^{s-4-\sigma(k,n)}},
\end{split}
\end{equation}
where $\sigma(k,n)$ is defined in \eqref{2.2}; 
we have assumed $s\geq[{{n-1}\over4}]+1$ and 
$s\geq[{{n-1}\over4}]+4$ for $M_0(T)$ and $M_1(T)$, respectively. 
The summations in $M_0(T)$ and $M_1(T)$ are taken over all 
$k\geq 0$ with $\sigma(k,n)\leq s-1$ and $\sigma(k,n)\leq s-4$, 
respectively. 
$E(T)$ is a time-weighted energy norm and $D(T)$ is the 
associated dissipation norm, while $M_0(T)$ and $M_1(T)$ 
are corresponding to the optimal decay estimates for $u$ and $u_t$, 
respectively. 

Now, the result on our a priori estimates is stated as follows. 

\bp[A priori estimates] \label{AP} 
Suppose that the conditions {\rm [A1]} and {\rm [A2]} are satisfied. 
Let $n\geq 2$ and $s\geq s(n)$, and assume that 
$u_0\in H^{s+1}(\R^n)\cap L^1(\R^n)$ and 
$u_1\in H^s(\R^n)\cap L^1(\R^n)$. 
Let $T>0$ and let $u(x,t)$ be the corresponding solution to the 
problem \eqref{1a}, \eqref{IC} satisfying 
$$
u\in C^0([0,T]; H^{s+1}(\R^n))\cap C^1([0,T]; H^s(\R^n)) 
$$
and the a priori bound \eqref{3a} below. 
Then there is a positive constant $\delta_1$ independent of $T$ 
such that if $E_1\leq\delta_1$, then the solution verifies the 
time-weighted energy estimate 
\begin{equation}\label{energy}
E(T)^2+D(T)^2\leq CE_1^2
\end{equation}
and the following optimal decay estimate: 
\begin{equation}\label{decay}
M_0(T)+M_1(T)\leq CE_1.
\end{equation}
Here $E_1$ is given in Theorem \ref{21}. 
\ep

\noindent
{\bf Remark.}\ 
In order to obtain the above a priori estimates, we employ a 
time-weighted energy method. To close our energy estimates, 
we make use of the following decay estimates for the $L^{\infty}$ 
norm of the derivatives 
$\p^2_xu_t$, $\p^3_xu$ and $\p^2_xu$: 
\begin{eqnarray*}
&\|\p^2_xu_t(t)\|_{L^\infty}
+\|\p^3_xu(t)\|_{L^{\infty}}\leq
CE_1(1+t)^{-\gamma},\ \ \gamma>1, \\[2mm]
&\|\p^2_xu(t)\|_{L^{\infty}}\leq CE_1(1+t)^{-d}, \ \ 
d>{n\over8}+{1\over2}.
\end{eqnarray*}
For the details, see Sections 3, 4 and 5. 

\medskip
Our next result is concerning the asymptotic profile of the global 
solution obtained in Theorem \ref{21}. 
First we show that the solution to the problem \eqref{1a}, \eqref{IC} 
can be approximated by the solution to the corresponding linear 
problem \eqref{1cc}, \eqref{IC}. Then we prove that the solution to 
this linear problem can be further approximated by the profile 
$MG_0(x,t+1)$, where 
$M=\int_{{\mathbb R}^n}(u_0+u_1)(x)dx$ denotes the "mass" and 
\begin{equation}\label{2.7}
G_0(x,t)=\mathcal{F}^{-1}[e^{-\gamma(\omega)|\xi|^4t}](x)
\end{equation}
is the fundamental solution to the fourth-order linear parabolic 
equation \eqref{limit-eq}, where $\gamma(\omega)$ is defined in [A2]. 
Notice that $G_0(x,t)$ is the self-similar solution of 
\eqref{limit-eq}. In fact, we have 
$G_0(x,t)=t^{-{n\over4}}\phi_0(x/t^{{1\over4}})$, where 
\begin{equation}\label{2.9}
\phi_0(x)=G_0(x,1)={\mathcal F}^{-1}[e^{-\gamma(\omega)|\xi|^4}](x).
\end{equation}
Thus we conclude that $MG_0(x,t+1)$ is an asymptotic profile of 
the solution to our problem \eqref{1a}, \eqref{IC}. 
This result on the asymptotic profile is stated as follows. 

\bt[Asymptotic profile]\label{22} 
Suppose that {\rm [A1]} and {\rm [A2]} are 
satisfied. Let $n\geq 2$ and $s\geq s(n)$. 
Assume that $u_0\in H^{s+1}(\R^n)\cap L^1_1(\R^n)$ and 
$u_1\in H^s(\R^n)\cap L^1_1(\R^n)$, and put 
$$
E_2:=\|u_0\|_{H^{s+1}}+\|u_1\|_{H^s}+\|(u_0,u_1)\|_{L^1_1}.
$$ 
Then the global solution $u(x,t)$ to the problem \eqref{1a}, 
\eqref{IC}, which is constructed in Theorem \ref{21}, is asymptotic 
to the profile $MG_0(x,t+1)$ as $t\to\infty$ in the following sense: 
\begin{equation}\label{asymp-1}
\|\p^k_x\{u(t)-MG_0(t+1)\}\|_{H^{s-2-\sigma(k,n)}}\leq
CE_2(1+t)^{-{n\over8}-{{k+1}\over4}},
\end{equation}
\begin{equation}\label{asymp-2}
\|\p^k_x\p_t\{u(t)-MG_0(t+1)\}\|_{H^{s-6-\sigma(k,n)}}\leq
CE_2(1+t)^{-{n\over8}-{{k+5}\over4}}
\end{equation}
for $k\geq 0$, where $\sigma(k,n)\leq s-2$ in \eqref{asymp-1} and 
$\sigma(k,n)\leq s-6$ in \eqref{asymp-2}. Here $M$ is a constant 
given by $M=\int_{{\mathbb R}^n}(u_0+u_1)(x)dx$ and $G_0(x,t)$ is 
the fundamental solution of \eqref{limit-eq} given in \eqref{2.7}. 
\et


\section{Time-weighted energy method}

In this section, we introduce a time-weighted energy method for our 
nonlinear problem \eqref{1a}, \eqref{IC} and explain why the standard 
energy method does not work well for our problem. 
First we give a lemma which will be used in the next energy estimates. 

\bl\label{31} 
Let $n\geq 1$, $1\leq p,q,r\leq\infty$ and 
${1\over p}={1\over q}+{1\over r}$. Then the following estimates hold: 
\be\label{d1} 
\|\p^k_x(uv)\|_{L^p}\leq
C(\|u\|_{L^q}\|\p^k_xv\|_{L^r}+\|v\|_{L^q}\|\p^k_xu\|_{L^r})
\ee
for $k\geq 0$, and 
\be\label{d2} 
\|[\p^k_x,u]\p_xv\|_{L^p}\leq
C(\|\p_xu\|_{L^q}\|\p^k_xv\|_{L^r}+\|\p_xv\|_{L^q}\|\p^k_xu\|_{L^r})
\ee
for $k\geq 1$, where $[A,B]=AB-BA$ denotes the commutator. 
\el

\begin{proof}
These estimates can be found in a literature but we give here a proof. 
To prove \eqref{d1}, it is enough to show that, for $k_1\geq1$, 
$k_2\geq1$ and $k_1+k_2=k$, the following estimate holds: 
$$
\|\p^{k_1}_xu\,\p^{k_2}_xv\|_{L^p}\leq
C(\|u\|_{L^q}\|\p^k_xv\|_{L^r}+\|v\|_{L^q}\|\p^k_xu\|_{L^r}).
$$
Let $\theta_j={k_j\over k}$, $j=1,2$, and define $p_j$, $j=1,2$, by 
$$
{1\over p_j}-{k_j\over n}
=(1-\theta_j){1\over q}+\theta_j({1\over r}-{k\over n}).
$$ 
Since $\theta_1+\theta_2=1$, we have 
${1\over p}={1\over p_1}+{1\over p_2}$. 
By using the H\"{o}lder inequality and the Gagliardo-Nirenberg 
inequality, we have 
\begin{equation*}
\begin{split}
\|\p^{k_1}_xu\,\p^{k_2}_xv\|_{L^p}\leq &\ 
\|\p^{k_1}_xu\|_{L^{p_1}}\|\p^{k_2}_xv\|_{L^{p_2}}\\[1mm]
\leq &\ 
C(\|u\|^{1-\theta_1}_{L^{q}}\|\p^{k}_xu\|^{\theta_1}_{L^{r}})
(\|v\|^{1-\theta_2}_{L^{q}}\|\p^{k}_xv\|^{\theta_2}_{L^{r}})\\[1mm]
\leq &\ 
C(\|u\|_{L^{q}}\|\p^{k}_xv\|_{L^{r}})^{\theta_2}
(\|v\|_{L^{q}}\|\p^{k}_xu\|_{L^{r}})^{\theta_1}\\[1mm]
\leq &\ 
C(\|u\|_{L^q}\|\p^k_xv\|_{L^r}+\|v\|_{L^q}\|\p^k_xu\|_{L^r}).
\end{split}
\end{equation*}
In the last inequality, we have used the Young inequality. Thus 
\eqref{d1} is proved. 

Next we show \eqref{d2}. We observe that 
$[\p^k_x,u]\p_xv$ consists of terms of the form 
$\p^{k_1}_x(\p_xu)\p^{k_2}_x(\p_xv)$ with $k_1$ and $k_2$ satisfying 
$k_1\geq0$, $k_2\geq0$ and $k_1+k_2=k-1$. 
Therefore, applying \eqref{d1}, we have 
$$
\|\p^{k_1}_x(\p_xu)\p^{k_2}_x(\p_xv)\|_{L^p}\leq
C(\|\p_xu\|_{L^q}\|\p^k_xv\|_{L^r}+\|\p_xv\|_{L^q}\|\p^k_xu\|_{L^r}),
$$
which gives \eqref{d2}. This completes the proof of Lemma \ref{31}.
\end{proof}

\medskip

{\bf Energy estimates:}\ \ 
Now, let $T>0$ and consider solutions to the problem \eqref{1a}, 
\eqref{IC}, which are defined on the time interval $[0,T]$ and 
verify the regularity mentioned in Proposition \ref{AP}. 
We derive energy estimates for the solutions under the following 
a priori assumption: 
\be\label{3a}
\sup_{0\leq t\leq T}\|\p^2_xu(t)\|_{L^{\infty}}\leq  \bar{\d},
\ee
where $\bar{\d}>0$ is a given small number not depending on $T$. 

First, we multiply the equation \eqref{1a} by $u_t$. 
After straightforward computations, we have the energy equality 
\begin{equation*}
\begin{split}
&\frac{1}{2}\{u_t^2+|\nabla u_t|^2+\phi(\partial_x^2u)\}_t
  +u_t^2-\nabla\cdot(u_t\nabla u_{tt}) \\
&+\sum_{ij}\{b^{ij}(\partial_x^2u)_{x_i}u_t\}_{x_j}
  -\sum_{ij}\{b^{ij}(\partial_x^2u)u_{tx_j}\}_{x_i}=0,
\end{split}
\end{equation*}
where $\phi(V)$ is the potential function in [A1]. 
We integrate this equality in $x\in{\mathbb R}^n$, obtaining 
\begin{equation}\label{3.4}
\frac{1}{2}\frac{d}{dt}\big\{
\|u_t\|_{H^1}^2+\int_{\mathbb{R}^n}\!\phi(\partial_x^2u)dx\big\}
+\|u_t\|_{L^2}^2=0.
\end{equation}
Here the Taylor expansion, using \eqref{1.6}, shows that 
$$
\phi(\partial_x^2u)=\frac{1}{2}\sum_{ij}\sum_{\alpha\beta}
b^{ij}_{\alpha\beta}(O)u_{x_ix_j}u_{x_\alpha x_\beta}
+O(|\partial_x^2u|^3), 
$$
which together with \eqref{elliptic} and \eqref{3a} gives 
\begin{equation}\label{3.5}
\int_{\mathbb{R}^n}\!\phi(\partial_x^2u)dx
\geq c\|\partial_x^2u\|_{L^2}^2.
\end{equation}

Next we derive a similar energy equality for derivatives. 
Notice that 
$b^{ij}(\p^2_xu)_{x_ix_j}=\sum_{\alpha\beta}
\{b^{ij}_{\a\b}(\p^2_xu)u_{x_{\a}x_{\b}x_i}\}_{x_j}$ in \eqref{1a}. 
Then, applying $\p^l_x$ to \eqref{1a}, we have 
\begin{equation}\label{3b}
\begin{split}
&\p^l_xu_{tt}-\p^l_x\Delta u_{tt}
  +\sum_{ij}\sum_{\alpha\beta}\{b^{ij}_{\a\b}(\p^2_xu)
  \p^l_xu_{x_{\a}x_{\b}x_i} \\
&+[\p^l_x,b^{ij}_{\a\b}(\p^2_xu)]u_{x_{\a}x_{\b}x_i}\}_{x_j}
  +\p^l_xu_t=0.
\end{split}
\end{equation}
We multiply this equation by $\p^l_xu_t$. After direct computations, 
we obtain 
\begin{equation*}\label{}
\begin{split}
&\frac{1}{2}\big\{|\p^l_xu_t|^2+|\nabla\p^l_xu_t|^2
  +\sum_{ij}\sum_{\alpha\beta}b^{ij}_{\a\b}(\p^2_xu)
  \p^l_xu_{x_ix_j}\p^l_xu_{x_{\a}x_{\b}}\big\}_t \\
&+|\p^l_xu_t|^2-\nabla\cdot(\p^l_xu_t\nabla\p^l_xu_{tt}) 
+\sum_{ij}\{\cdots\}_{x_j}-\sum_{ij}\{\cdots\}_{x_i}=r^{(l)},
\end{split}
\end{equation*}
where 
\begin{equation*}
\begin{split}
r^{(l)}
&=\sum_{ij}\sum_{\a\b}\big\{
  \frac{1}{2}\,b^{ij}_{\a\b}(\p^2_xu)_t\,
  \p^l_xu_{x_{\a}x_{\b}}\p^l_xu_{x_ix_j} \\
&-b^{ij}_{\a\b}(\p^2_xu)_{x_i}\p^l_xu_{x_{\a}x_{\b}}\p^l_xu_{tx_j} 
  +[\p^l_x,b^{ij}_{\a\b}(\p^2_xu)]
  u_{x_{\a}x_{\b}x_i}\p^l_xu_{tx_j}\big\}.
\end{split}
\end{equation*}
We integrate the above equality in $x\in{\mathbb R}^n$, obtaining 
\be\label{3c}
\frac{1}{2}\frac{d}{dt}
  \{\|\p^l_xu_t\|^2_{H^1}+\|\p^{l}_xu\|^2_{L^2(b)}\}
  +\|\p^l_xu_t\|_{L^2}^2\leq R^{(l)}
\ee
where we put $R^{(l)}=\int_{\mathbb{R}^n}|r^{(l)}|\,dx$ and 
\begin{equation*}
\|v\|_{L^2(b)}^2=\sum_{ij}\sum_{\alpha\beta}\int_{\mathbb{R}^n}
  b^{ij}_{\a\b}(\p^2_xu)v_{x_ix_j}v_{x_{\a}x_{\b}}dx.
\end{equation*}
Here, using \eqref{elliptic} and \eqref{3a}, we see that 
\begin{equation}\label{3.8}
\|\p^{l}_xu\|^2_{L^2(b)}\geq c\|\p^{l+2}_xu\|^2_{L^2}.
\end{equation}
Also, applying Lemma \ref{31} together with \eqref{3a}, we find that 
$\|[\p^l_x,b(\p^2_xu)]\p^3_xu\|_{L^2} \leq C\|\p^3_xu\|_{L^\infty}
\|\p^{l+2}_xu\|_{L^2}$. Therefore the term $R^{(l)}$ can be 
estimated as 
\begin{equation}\label{3.9}
R^{(l)}\leq C(\|\p^2_xu_t\|_{L^\infty}+\|\p^3_xu\|_{L^\infty})
(\|\p^{l+1}_xu_t\|_{L^2}^2+\|\p^{l+2}_xu\|_{L^2}^2).
\end{equation}
We can use \eqref{3c} even for $l=0$ instead of the energy 
equality \eqref{3.4}. 

Finally, we multiply \eqref{1a} by $u$. This yields 
\begin{equation*}
\begin{split}
&\frac{1}{2}\{u^2+2u_t(u-\Delta u)\}_t
  +\sum_{ij}b^{ij}(\p^2_xu)u_{x_ix_j}-(u_t^2+|\nabla u_t|^2) \\
& -\nabla\cdot\{u\nabla u_{tt}-(u_t\nabla u)_t\}
  +\sum_{ij}\{\cdots\}_{x_j}-\sum_{ij}\{\cdots\}_{x_i}=0.
\end{split}
\end{equation*}
We integrate this equality in $x\in{\mathbb R}^n$. 
Since $b^{ij}(\p^2_xu)=\sum_{\a\b}b^{ij}_{\a\b}(O)u_{x_\a x_\b}
+O(|\p^2_xu|^2)$, we obtain 
\begin{equation}\label{3.10}
\frac{1}{2}\{\|u\|_{L^2}^2
+2\langle u_t,u-\Delta u\rangle_{L^2}\}_t
+c\|\p^2_xu\|_{L^2}-\|u_t\|_{H^1}^2\leq 0,
\end{equation}
where we have used \eqref{elliptic} and \eqref{3a}. 
To get a similar estimate for derivatives, we multiply \eqref{3b} 
by $\p^l_xu$. After direct computations, we have 
\begin{equation*}
\begin{split}
&\frac{1}{2}\{|\p^l_xu|^2+2\p^l_xu_t(\p^l_xu-\Delta\p^l_xu)\}_t
  +\sum_{ij}\sum_{\a\b}b^{ij}_{\a\b}(\p^2_xu)
  \p^l_xu_{x_ix_j}\p^l_xu_{x_\a x_\b} \\
&-(|\p^l_xu_t|^2+|\nabla\p^l_xu_t|^2) 
  +\sum_{ij}\{\cdots\}_{x_j}-\sum_{ij}\{\cdots\}_{x_i}
  =\tilde{r}^{(l)},
\end{split}
\end{equation*}
where 
\begin{equation*}
\tilde{r}^{(l)}
=\sum_{ij}\sum_{\a\b}\big\{
  -b^{ij}_{\a\b}(\p^2_xu)_{x_i}\p^l_xu_{x_{\a}x_{\b}}\p^l_xu_{x_j} 
  +[\p^l_x,b^{ij}_{\a\b}(\p^2_xu)]
  u_{x_{\a}x_{\b}x_i}\p^l_xu_{x_j}\big\}.
\end{equation*}
Integrating the above equality in $x\in{\mathbb R}^n$ and using 
\eqref{elliptic} and \eqref{3a}, we obtain 

\begin{equation}\label{3d}
\begin{split}
&\frac{1}{2}\frac{d}{dt}\{\|\p^l_xu\|^2_{L^2}
  +2\langle\p^{l}_xu_t,\p^l_xu-\Delta\p^{l}_xu\rangle_{L^2}\} \\[1mm]
&+c\|\p^{l+2}_xu\|^2_{L^2}-\|\p^l_xu_t\|_{H^1}^2
  \leq\tilde{R}^{(l)},
\end{split}
\end{equation}
where $\tilde{R}^{(l)}=\int_{\mathbb{R}^n}|\tilde{r}^{(l)}|\,dx$. 
The term $\tilde{R}^{(l)}$ is estimated as 
\begin{equation}\label{3.12}
\tilde{R}^{(l)}\leq C\|\p^3_xu\|_{L^\infty}
\|\p^{l+2}_xu\|_{L^2}\|\p^{l+1}_xu\|_{L^2}.
\end{equation}
Notice that \eqref{3.10} coincides with \eqref{3d} with $l=0$ 
if we set $\tilde{R}^{(0)}=0$. 

\smallskip
We explain why the standard energy method does not work well 
for our problem. To this end, we integrate \eqref{3c} with respect 
to $t$ and add the resulting inequality for $l$ with 
$0\leq l\leq s-1$. This yields 
\begin{equation}\label{3e}
\begin{split}
&\|u_t(t)\|^2_{H^s}+\|\p^2_xu(t)\|^2_{H^{s-1}}
  +\int^t_0\|u_t(\tau)\|^2_{H^{s-1}}d\tau \\
&\leq C(\|u_0\|^2_{H^{s+1}}+\|u_1\|^2_{H^s})+CR_1(t), 
\end{split}
\end{equation}
where $R_1(t)=\int_0^t\sum_{l=0}^{s-1}R^{(l)}(\tau)d\tau$. 
Also, we integrate \eqref{3.10} and \eqref{3d} with respect to $t$ 
and add the resulting inequalities for $l$ with $1\leq l\leq s-2$; 
we use \eqref{3d} only for $s\geq 3$. This gives 
\begin{equation}\label{3f}
\begin{split}
&\|u(t)\|^2_{H^{s-2}}+\int^t_0\|\p^2_xu(\tau)\|^2_{H^{s-2}}d\tau \\
&\leq C(\|u_0\|^2_{H^{s}}+\|u_1\|^2_{H^{s-2}})
  +C(\|u_t(t)\|^2_{H^{s-2}}+\|\p^2_xu(t)\|^2_{H^{s-2}}) \\
&+C\!\int^t_0\|u_t(\tau)\|^2_{H^{s-1}}d\tau+CR_2(t),
\end{split}
\end{equation}
where $R_2(t)=\int_0^t\sum_{l=0}^{s-2}\tilde{R}^{(l)}(\tau)d\tau$, 
in which we can regard $\tilde{R}^{(0)}=0$. 
Combining \eqref{3e} and \eqref{3f}, we arrive at the final energy 
inequality 
\begin{equation}\label{3g}
\begin{split}
&\|u_t(t)\|^2_{H^s}+\|u(t)\|^2_{H^{s+1}}
  +\int^t_0(\|u_t(\tau)\|^2_{H^{s-1}}
  +\|\p^2_xu(\tau)\|^2_{H^{s-2}})d\tau \\
&\leq C(\|u_0\|^2_{H^{s+1}}+\|u_1\|^2_{H^s})+CR(t), 
\end{split}
\end{equation}
where $R(t)=R_1(t)+R_2(t)$. We note that 
$$
R(t)\leq C\!\int^t_0\|(\p^2_xu_t,\p^3_xu)(\tau)\|_{L^{\infty}}
(\|u_t(\tau)\|^2_{H^s}+\|\p^2_xu(\tau)\|^2_{H^{s-1}})d\tau.
$$

To complete our energy method, we need to control the term $R(t)$ 
which comes from the nonlinearity of our equation \eqref{1a}. 
Usually, this can be done by using the dissipative term, namely, 
the third term on the left hand side of \eqref{3g}. 
Following to this strategy, we estimate the term $R(t)$ as 
$$
R(t)\leq C\sup\limits_{0\leq\tau\leq t}
\|(\p^2_xu_t,\p^3_xu)(\tau)\|_{L^{\infty}}\!
\int^t_0(\|u_t(\tau)\|^2_{H^{s}}
+\|\p^2_xu(\tau)\|^2_{H^{s-1}})d\tau.
$$
In our case, however, the dissipative term does not contain 
the highest-order term $\int^t_0(\|\p^s_xu_t(\tau)\|^2_{L^2}+
\|\p^{s+1}_xu(\tau)\|^2_{L^2})d\tau$ because of the loss of 
regularity, and therefore it could not control the nonlinearity 
$R(t)$. Consequently, the standard energy method does not work well 
for our problem. 

\medskip

{\bf Time-weighted energy estimates:}\ \ 
To resolve the above difficulty caused by the regularity-loss 
property, we try to estimate the nonlinearity $R(t)$ as 
$$
R(t)\leq C\sup\limits_{0\leq\tau\leq t}
(\|u_t(\tau)\|^2_{H^{s}}+\|\p^2_xu(\tau)\|^2_{H^{s-1}})\!
\int^t_0\|(\p^2_xu_t,\p^3_xu)(\tau)\|_{L^{\infty}}d\tau.
$$
This requires the integrability of 
$\|(\p^2_xu_t,\p^3_xu)(t)\|_{L^{\infty}}$ over $t\geq0$. 
To ensure this integrability, we need to show a suitable decay 
estimate for $\|(\p^2_xu_t,\p^3_xu)(t)\|_{L^{\infty}}$ and this 
will be done by employing the time-weighted energy method 
combined with the optimal $L^2$ decay estimates for lower-order 
derivatives of solutions. 

Now we estimate the time-weighted energy norm $E(T)$ and the 
associated dissipation norm $D(T)$ defined in \eqref{norm-1} 
by applying the time-weighted energy method mentioned above. 
We make use of the following integral norm $L(T)$: 
\begin{equation}
L(T):=\int^T_0\|(\p^2_xu_t,\p^3_xu)(\tau)\|_{L^{\infty}}d\tau.
\end{equation}
The result is stated as follows. 

\bp \label{p41} 
Suppose that the conditions {\rm [A1]} and {\rm [A2]} are satisfied. 
Let $n\geq 1$ and $s\geq 2$. Assume that $u_0\in H^{s+1}(\R^n)$ and 
$u_1\in H^s(\R^n)$, and put $E_0:=\|u_0\|_{H^{s+1}}+\|u_1\|_{H^s}$. 
Let $u(x,t)$ be the corresponding solution to the problem \eqref{1a}, 
\eqref{IC} which is defined on $[0,T]$ and verifies \eqref{3a}. 
Then we have the following estimate: 
$$
E(T)^2+D(T)^2\leq CE_0^2+CL(T)E(T)^2.
$$ 
\ep

\begin{proof}
For the proof, it is enough to show the following estimates 
for any $t\in[0,T]$: 
\begin{equation}\label{4a}
\|u(t)\|^2_{H^{s+1}}+\|u_t(t)\|^2_{H^s}
+\int^t_0\|u_t(\tau)\|^2_{H^{s-1}}d\tau\leq
CE_0^2+CL(T)E(T)^2,
\end{equation}
\begin{equation}\label{4b}
\int^t_0(1+\tau)^{j-1}\|\p^{2j}_xu(\tau)\|^2_{H^{s-3j+1}}d\tau\leq
CE_0^2+CL(T)E(T)^2,
\end{equation}
\begin{equation}\label{4c}
\begin{split}\\[-3mm]
&(1+t)^j(\|\p^{2j-2}_xu_t(t)\|^2_{H^{s-3j+2}}
  +\|\p^{2j}_xu(t)\|^2_{H^{s-3j+1}})\\
&+\int^t_0(1+\tau)^j\|\p^{2j-2}_xu_t(\tau)\|^2_{H^{s-3j+1}}d\tau
  \leq CE_0^2+CL(T)E(T)^2,
\end{split}
\end{equation}
where $1\leq j\leq [{{s+1}\over3}]$. 

We use \eqref{3c} and \eqref{3d}. 
Let $\lambda\geq 0$. We multiply \eqref{3c} by $(1+t)^\lambda$ 
and integrate with respect to $t$. This yields 
\begin{equation}\label{3x}
\begin{split}
&(1+t)^\lambda
  (\|\p^l_xu_t(t)\|^2_{H^1}+\|\p^{l+2}_xu(t)\|^2_{L^2})
  +\int^t_0(1+\tau)^\lambda
  \|\p^l_xu_t(\tau)\|^2_{L^2}d\tau \\
&\leq CE_0^2
  +\lambda C\!\int^t_0(1+\tau)^{\lambda-1}
  (\|\p^l_xu_t(\tau)\|^2_{H^1}+\|\p^{l+2}_xu(\tau)\|^2_{L^2})d\tau \\
&\qquad\qquad\qquad
  +C\!\int_0^t(1+\tau)^\lambda R^{(l)}(\tau)d\tau
\end{split}
\end{equation}
for $0\leq l\leq s-1$. Here we can regard $R^{(0)}=0$ 
if we use \eqref{3.4} in place of \eqref{3c} with $l=0$. 
Also, we multiply \eqref{3d} by $(1+t)^\lambda$ 
and integrate with respect to $t$. This yields 
\begin{equation}\label{3y}
\begin{split}
&(1+t)^\lambda\|\p^l_xu(t)\|^2_{L^2}
  +\int^t_0(1+\tau)^\lambda\|\p^{l+2}_xu(\tau)\|^2_{L^2}d\tau \\
&\leq CE_0^2
  +\lambda C\!\int^t_0(1+\tau)^{\lambda-1}
  (\|\p^l_xu_t(\tau)\|^2_{L^2}+\|\p^{l}_xu(\tau)\|^2_{H^2})d\tau \\
&+C(1+t)^\lambda
  (\|\p^l_xu_t(t)\|^2_{L^2}+\|\p^{l+2}_xu(t)\|^2_{L^2})\\
&+C\!\int_0^t(1+\tau)^\lambda\|\p^l_xu_t(\tau)\|^2_{H^1}d\tau
  +C\!\int_0^t(1+\tau)^\lambda\tilde{R}^{(l)}(\tau)d\tau
\end{split}
\end{equation}
for $0\leq l\leq s-2$. Here we can regard $\tilde{R}^{(0)}=0$ 
if we use \eqref{3.10} in place of \eqref{3d} with $l=0$. 

First we show \eqref{4a} and \eqref{4b} with $j=1$. 
To this end, we put $\lambda=0$ in \eqref{3x} and add for $l$ 
with $0\leq l\leq s-1$. Since 
$\int_0^t\sum_{l=0}^{s-1}R^{(l)}(\tau)d\tau\leq CL(T)E(T)^2$, 
we have 
\begin{equation}\label{3.e}
\begin{split}
&\|u_t(t)\|^2_{H^s}+\|\p^2_xu(t)\|^2_{H^{s-1}}
  +\int^t_0\|u_t(\tau)\|^2_{H^{s-1}}d\tau \\
&\leq CE_0^2+CL(T)E(T)^2,
\end{split}
\end{equation}
which corresponds to \eqref{3e}. Also, we put $\lambda=0$ in 
\eqref{3y} and add for $l$ with $0\leq l\leq s-2$. Since 
$\int_0^t\sum_{l=0}^{s-2}\tilde{R}^{(l)}(\tau)d\tau\leq CL(T)E(T)^2$, 
we obtain 
\begin{equation}\label{3.f}
\begin{split}
&\|u(t)\|^2_{H^{s-2}}
  +\int^t_0\|\p^2_xu(\tau)\|^2_{H^{s-2}}d\tau \\[-1mm]
&\leq CE_0^2
  +C(\|u_t(t)\|^2_{H^{s-2}}+\|\p^2_xu(t)\|^2_{H^{s-2}}) 
  +C\!\int^t_0\|u_t(\tau)\|^2_{H^{s-1}}d\tau \\[1mm]
&\qquad
  +CL(T)E(T)^2 \leq CE_0^2+CL(T)E(T)^2,
\end{split}
\end{equation}
where we have used \eqref{3.e} in the last estimate. 
This is corresponding to \eqref{3f}. 
The desired estimates \eqref{4a} and \eqref{4b} with $j=1$ follow 
from \eqref{3.e} and \eqref{3.f}. 

Next we show \eqref{4c} with $j=1$. To this end, we put 
$\lambda=1$ in \eqref{3x} and add for $l$ with $0\leq l\leq s-2$. 
This yields 
\begin{equation*}
\begin{split}
&(1+t)(\|u_t(t)\|^2_{H^{s-1}}+\|\p^2_xu(t)\|^2_{H^{s-2}})
  +\int^t_0(1+\tau)\|u_t(\tau)\|^2_{H^{s-2}}d\tau \\
&\leq CE_0^2
  +C\!\int_0^t(\|u_t(\tau)\|^2_{H^{s-1}}
  +\|\p^2_xu(\tau)\|^2_{H^{s-2}})d\tau \\
&\qquad
  +C\!\int_0^t(1+\tau)\sum_{l=0}^{s-2}R^{(l)}(\tau)d\tau 
  \leq CE_0^2+CL(T)E(T)^2,
\end{split}
\end{equation*}
where we have used \eqref{3.e}, \eqref{3.f} and the fact that 
the last integral involving $R^{(l)}$ can be estimated by 
\begin{equation*}
\begin{split}
&C\!\sup_{0\leq\tau\leq t}
  (1+\tau)(\|u_t(\tau)\|^2_{H^{s-1}}+\|\p^2_xu(\tau)\|^2_{H^{s-2}}) 
  \int_0^t\|(\p^2_xu_t,\p^3_xu)(\tau)\|_{L^\infty}d\tau \\
&\leq CL(T)E(T)^2.
\end{split}
\end{equation*}
Thus we have shown \eqref{4c} with $j=1$. 

Now we show \eqref{4b} and \eqref{4c} by induction with respect 
to $j$. For this purpose, let $s\geq 5$ and 
$1\leq k\leq[{{s+1}\over3}]-1=[{{s-2}\over3}]$, and assume 
that (\ref{4b}) and (\ref{4c}) hold true for $j=k$. 
Then we show (\ref{4b}) and (\ref{4c}) for $j=k+1$. 
To this end, we first put $\lambda=k$ in \eqref{3y} and add 
for $l$ with $2k\leq l \leq s-k-2$. This yields 
\begin{equation*}
\begin{split}
&(1+t)^k\|\p^{2k}_xu(t)\|^2_{H^{s-3k-2}}
  +\int^t_0(1+\tau)^k\|\p^{2k+2}_xu(\tau)\|^2_{H^{s-3k-2}}d\tau \\
&\leq CE_0^2
  +C\!\int^t_0(1+\tau)^{k-1}(\|\p^{2k}_xu_t(\tau)\|^2_{H^{s-3k-2}}
  +\|\p^{2k}_xu(\tau)\|^2_{H^{s-3k}})d\tau \\
&+C(1+t)^k(\|\p^{2k}_xu_t(t)\|^2_{H^{s-3k-2}}
  +\|\p^{2k+2}_xu(t)\|^2_{H^{s-3k-2}}) \\
&+C\!\int^t_0(1+\tau)^k\|\p^{2k}_xu_t(\tau)\|^2_{H^{s-3k-1}}d\tau 
  +C\!\int^t_0(1+\tau)^k\sum_{l=2k}^{s-k-2}
  \tilde{R}^{(l)}(\tau)\,d\tau.
\end{split}
\end{equation*}
Here the last integral can be estimated by 
\begin{equation*}
C\!\sup_{0\leq\tau\leq t}
  (1+\tau)^k\|\p^{2k}_xu(\tau)\|^2_{H^{s-3k}}\! 
  \int_0^t\|\p^3_xu(\tau)\|_{L^\infty}d\tau 
\leq CL(T)E(T)^2.
\end{equation*}
Also, the term $C\!\int^t_0(1+\tau)^{k-1}
\|\p^{2k}_xu(\tau)\|^2_{H^{s-3k}}d\tau$ is estimated by using 
\eqref{4b} with $j=k$, while the other terms on the right hand 
side of the above inequality are estimated by making use of 
\eqref{4c} with $j=k$. Consequently, we obtain 
\begin{equation*}
\begin{split}
&(1+t)^k\|\p^{2k}_xu(t)\|^2_{H^{s-3k-2}}
  +\int^t_0(1+\tau)^k\|\p^{2k+2}_xu(\tau)\|^2_{H^{s-3k-2}}d\tau \\
&\leq CE_0^2+CL(T)E(T)^2,
\end{split}
\end{equation*}
which implies the desired estimate \eqref{4b} for $j=k+1$. 

Finally, we show \eqref{4c} for $j=k+1$. 
Let $\lambda=k+1$ in \eqref{3x} and add for $l$ with 
$2k\leq l\leq s-k-2$. This yields 
\begin{equation*}
\begin{split}
&(1+t)^{k+1}(\|\p^{2k}_xu_t(t)\|^2_{H^{s-3k-1}}
  +\|\p^{2k+2}_xu(t)\|^2_{H^{s-3k-2}}) \\
&\qquad+\int^t_0(1+\tau)^{k+1}
  \|\p^{2k}_xu_t(\tau)\|^2_{H^{s-3k-2}}d\tau \\
&\leq CE_0^2
  +C\!\int^t_0(1+\tau)^k
  (\|\p^{2k}_xu_t(\tau)\|^2_{H^{s-3k-1}}
  +\|\p^{2k+2}_xu(\tau)\|^2_{H^{s-3k-2}})d\tau \\
&\qquad
  +C\!\int^t_0(1+\tau)^{k+1}
  \sum_{l=2k}^{s-k-2}R^{(l)}(\tau)d\tau.
\end{split}
\end{equation*}
Here the last integral can be estimated by 
\begin{equation*}
\begin{split}
&C\!\sup_{0\leq\tau\leq t}(1+\tau)^{k+1}
  (\|\p^{2k}u_t(\tau)\|^2_{H^{s-3k-1}}
  +\|\p^{2k+2}_xu(\tau)\|^2_{H^{s-3k-2}}) \\
&\qquad
  \cdot\int_0^t\|(\p^2_xu_t,\p^3_xu)(\tau)\|_{L^\infty}d\tau 
  \leq CL(T)E(T)^2.
\end{split}
\end{equation*}
Also, the term 
$C\!\int^t_0(1+\tau)^k\|\p^{2k}_xu_t(\tau)\|^2_{H^{s-3k-1}}d\tau$ 
is estimated by \eqref{4c} with $j=k$, while the term 
$C\!\int^t_0(1+\tau)^k\|\p^{2k+2}_xu(\tau)\|^2_{H^{s-3k-2}}d\tau$ 
can be estimated by using \eqref{4b} with $j=k+1$ which has already 
been proved. Thus we get 
\begin{equation*}
\begin{split}
&(1+t)^{k+1}(\|\p^{2k}_xu_t(t)\|^2_{H^{s-3k-1}}
  +\|\p^{2k+2}_xu(t)\|^2_{H^{s-3k-2}}) \\
&+\int^t_0(1+\tau)^{k+1}
  \|\p^{2k}_xu_t(\tau)\|^2_{H^{s-3k-2}}d\tau
  \leq CE_0^2+CL(T)E(T)^2,
\end{split}
\end{equation*}
which proves \eqref{4c} for $j=k+1$. This completes the proof of 
Proposition \ref{41}. 
\end{proof}


\section{Optimal decay estimates}

This section is devoted to the proof of the optimal decay estimates 
for solutions to the problem \eqref{1a}, \eqref{IC}, which are 
defined on the time interval $[0,T]$ and verify \eqref{3a}. 

Recalling \eqref{1.6}, we write $b^{ij}(V)$ as 
$b^{ij}(V)=\sum_{\alpha\beta}b^{ij}_{\alpha\beta}(O)V_{\alpha\beta}
+g^{ij}(V)$, where $g^{ij}(V)$ satisfies $g^{ij}(V)=O(|V|^2)$ 
for $|V|\to 0$ and the symmetric property \eqref{1.4}. 
The equation \eqref{1a} is then rewritten in the form 
\begin{equation}\label{4.1}
u_{tt}-\Delta u_{tt}
+\sum_{i,j=1}^n\sum_{\a,\b=1}^n b^{ij}_{\a\b}(O)u_{x_ix_jx_{\a}x_{\b}}
+\sum_{i,j=1}^n g^{ij}(\partial_x^2u)_{x_ix_j}+u_t=0.
\end{equation}
Consequently, by the Duhamel principle, we can express the solutions 
to the problem \eqref{1a}, \eqref{IC} as 
\begin{equation}\label{4.2}
\begin{split}
u(t)
&=G(t)\ast(u_0+u_1)+H(t)\ast u_0 \\
&-\int^t_0G(t-\tau)\ast(1-\Delta)^{-1}
  \p^2_xg(\p^2_xu)(\tau)d\tau,
\end{split}
\end{equation}
where $G(x,t)$ and $H(x,t)$ are the fundamental solutions to the 
linearized equation \eqref{1cc}. Here and in the following, we 
use the abbreviation 
$\p^2_xg(\p^2_xu)=\sum_{ij}g^{ij}(\p^2_xu)_{x_ix_j}$. 
The fundamental solutions $G(x,t)$ and $H(x,t)$ are given by 
the formulas 
$$
G(x,t)=\mathcal{F}^{-1}\Big[
{{e^{\lambda_+(\xi)t}-e^{\lambda_-(\xi)t}}
\over{\lambda_+(\xi)-\lambda_-(\xi)}}\Big](x),
$$
$$
H(x,t)=\mathcal{F}^{-1}\Big[
{{(1+\lambda_+(\xi))e^{\lambda_-(\xi)t}
-(1+\lambda_-(\xi))e^{\lambda_+(\xi)t}}
\over{\lambda_+(\xi)-\lambda_-(\xi)}}\Big](x),
$$
where $\lambda_\pm(\xi)$ are the eigenvalues given explicitly by 
$$
\lambda_{\pm}(\xi)
={{-1\pm\sqrt{1-4\gamma(\omega)|\xi|^4(1+|\xi|^2)}}\over{2(1+|\xi|^2)}}
$$
with $\gamma(\omega)$ being defined in [A2]. 
We decompose the solution formula \eqref{4.2} in the form 
$u(t)=\bar{u}(t)-F(u)(t)$, where $\bar{u}(x,t)$ and $F(u)(x,t)$ 
denote the linear and nonlinear parts, respectively: 
\begin{equation}\label{4.3}
\begin{split}
&\bar{u}(t)=G(t)\ast(u_0+u_1)+H(t)\ast u_0, \\
&F(u)(t)=\int^t_0G(t-\tau)\ast(1-\Delta)^{-1}
  \p^2_xg(\p^2_xu)(\tau)d\tau.
\end{split}
\end{equation}

We review the basic decay property for the equation \eqref{1c} 
which was studied in \cite{SK}. 
Notice that \eqref{1c} is a special case of our linearized equation 
\eqref{1cc} with $b^{ij}_{\a\b}(O)=\delta_{ij}\delta_{\a\b}$, 
in which we have $\gamma(\omega)=1$. 
Therefore the corresponding fundamental solutions to \eqref{1c} are 
given by the above formula with $\gamma(\omega)=1$. 
In particular, we have 
$$
\tilde{G}(x,t)=\mathcal{F}^{-1}\Big[
{{e^{\tilde{\lambda}_+(\xi)t}-e^{\tilde{\lambda}_-(\xi)t}}
\over{\tilde{\lambda}_+(\xi)-\tilde{\lambda}_-(\xi)}}\Big](x),
$$
where 
$$
\tilde{\lambda}_{\pm}(\xi)
={{-1\pm\sqrt{1-4|\xi|^4(1+|\xi|^2)}}\over{2(1+|\xi|^2)}}.
$$
The following decay result was proved in \cite{SK}. 

\bl[\cite{SK}]\label{41} 
Let $n\geq 1$ and $s\geq 0$, and assume that 
$\phi\in H^s(\R^n)\cap L^1(\R^n)$. 
Then the following decay estimates hold: 
\begin{equation*}
\|\p^k_x\tilde{G}(t)\ast\phi\|_{L^2}
\leq C(1+t)^{-{n\over8}-{k\over4}}\|\phi\|_{L^1}
+C(1+t)^{-{{l+1}\over2}}\|\p^{k+l}_x\phi\|_{L^2}
\end{equation*}
for integers $k$ and $l$ with $k\geq 0$, $l+1\geq0$ and 
$0\leq k+l\leq s$, and 
\begin{equation*}
\|\p^k_x\p_t\tilde{G}(t)\ast\phi\|_{L^2}
\leq C(1+t)^{-{n\over8}-{k\over4}-1}\|\phi\|_{L^1}
+C(1+t)^{-{{l}\over2}}\|\p^{k+l}_x\phi\|_{L^2}
\end{equation*}
for integers $k$ and $l$ with $k\geq0$, $l\geq0$ and $k+l\leq s$. 
\el

In view of the structural assumption $[A2]$, the quantity 
$\gamma(\omega)$ has a positive minimum over $S^{n-1}$. 
%
%
By virtue of this fact, we find that the proof of Lemma \ref{41} 
in \cite{SK} is valid also for our fundamental solution $G(x,t)$. 
Thus we conclude that our $G(x,t)$ satisfies the same decay estimates 
in Lemma \ref{41}. 
Moreover, as a simple modification of this lemma, we also have: 
\begin{equation}\label{d-1}
\begin{split}
&\|\p^{k}_xG(t)\ast(1-\Delta)^{-1}\phi\|_{L^2} \\
&\leq C(1+t)^{-{n\over8}-{{k}\over4}}\|\phi\|_{L^1}
  +C(1+t)^{-{{l+1}\over2}}\|\p^{k+l-2}_x\phi\|_{L^2}
\end{split}
\end{equation}
for $k\geq0$, $l+1\geq0$ and $2\leq k+l\leq s+2$, and 
\begin{equation}\label{d-2}
\begin{split}
&\|\p^{k}_x\p_tG(t)\ast(1-\Delta)^{-1}\phi\|_{L^2} \\
&\leq C(1+t)^{-{n\over8}-{{k}\over4}-1}\|\phi\|_{L^1}
  +C(1+t)^{-{{l}\over2}}\|\p^{k+l-2}_x\phi\|_{L^2}
\end{split}
\end{equation}
for $k\geq0$, $l\geq0$ and $2\leq k+l\leq s+2$.

As for the solution to our linearized problem, 
by similar calculations as in \cite{SK}, we obtain the following 
decay result. 

\bl\label{Su} 
Let $n\geq 1$ and assume that $u_0\in H^{s+1}(\R^n)\cap L^1(\R^n)$ 
and $u_1\in H^s(\R^n)\cap L^1(\R^n)$ for an integer $s$ specified 
below. Put 
$E_1=\|u_0\|_{H^{s+1}}+\|u_1\|_{H^{s}}+\|(u_0,u_1)\|_{L^1}$. 
Then the linear part $\bar{u}(x,t)$ in \eqref{4.3} satisfies 
the following decay estimates: 
\begin{equation}\label{dL1}
\|\p^k_x\bar{u}(t)\|_{H^{s-\sigma_1(k,n)}}\leq
CE_1(1+t)^{-{n\over8}-{k\over4}},
\end{equation}
\begin{equation}\label{dL2}
\|\p^k_x\bar{u}_t(t)\|_{H^{s-3-\sigma_1(k,n)}}\leq
CE_1(1+t)^{-{n\over8}-{k\over4}-1},
\end{equation}
where $s\geq[{{n-1}\over4}]$, $k\geq 0$ and $\sigma_1(k,n)\leq s$ 
in \eqref{dL1}, and 
$s\geq[{{n-1}\over4}]+3$, $k\geq 0$ and $\sigma_1(k,n)\leq s-3$ 
in \eqref{dL2}. Here $\sigma_1(k,n)$ is defined in \eqref{2.1}. 
\el 

Now we estimate the time-weighted norms $M_0(T)$ and $M_1(T)$ 
defined in \eqref{norm-2}, which are corresponding to the optimal 
$L^2$ decay estimates of solutions to the nonlinear problem. 
To control these norms, we will use the following time-weighted 
$L^\infty$ norm: 
\begin{equation}
N_d(T):=\sup\limits_{0\leq t\leq T}(1+t)^d
\|\p^2_xu(t)\|_{L^{\infty}},
\end{equation}
where $d>d(n):={n\over8}+{1\over2}$. 
By applying Lemmas \ref{41} and \ref{Su}, we have: 

\bp\label{p42} 
Suppose that {\rm [A1]} and {\rm [A2]} are 
satisfied. Let $n\geq 2$ and $s\geq[{{n-1}\over4}]+4$. 
Assume that $u_0\in H^{s+1}(\R^n)\cap L^1(\R^n)$ and 
$u_1\in H^s(\R^n)\cap L^1(\R^n)$, and put 
$E_1=\|u_0\|_{H^{s+1}}+\|u_1\|_{H^s}+\|(u_0,u_1)\|_{L^1}$. 
Let $u(x,t)$ be the corresponding solution to the problem \eqref{1a}, 
\eqref{IC} which is defined on $[0,T]$ and verifies \eqref{3a}. 
Then the following estimates hold: 
\begin{equation}\label{dM-0}
M_0(T)\leq CE_1+CM_0(T)^2+C(M_0(T)+N_d(T))E(T),
\end{equation}
\begin{equation}\label{dM-1}
M_1(T)\leq CE_1+CM_0(T)^2+C(M_0(T)+N_d(T))E(T).
\end{equation}
Here $d$ is assumed to satisfy $d>d(n)={n\over8}+{1\over2}$. 
\ep

\begin{proof} 
First we prove \eqref{dM-0}. Let $k\geq0$ and $h\geq0$ be integers. 
We apply $\p^{k+h}_x$ to $F(u)$ in \eqref{4.3} and take the $L^2$ 
norm, obtaining 
\begin{equation}\label{N-0}
\begin{split}
\|\p^{k+h}_xF(u)(t)\|_{L^2}
&\leq\int_0^t\|\p^{k+h+2}_xG(t-\tau)\ast(1-\Delta)^{-1}
  g(\p^2_xu)(\tau)\|_{L^2}d\tau \\
&=\int^{t\over2}_0+\int^t_{t\over2}=:I_1+I_2.
\end{split}
\end{equation}
For the term $I_1$, we apply \eqref{d-1} with $k$ replaced by 
$k+h+2$ and with $\phi=g(\p^2_xu)$ to get 
\be\label{I-1}
\begin{split}
I_1
&\leq C\!\int^{{t\over2}}_0(1+t-\tau)^{-{n\over8}-{{k+h+2}\over4}}
  \|g(\p^2_xu)(\tau)\|_{L^1}d\tau \\
&+C\!\int^{{t\over2}}_0(1+t-\tau)^{-{{l+1}\over2}}
  \|\p^{k+h+l}_xg(\p^2_xu)(\tau)\|_{L^2}d\tau
  =:I_{11}+I_{12}.
\end{split}
\ee
We see that $\|g(\p^2_xu)\|_{L^1}\leq C\|\p^2_xu\|_{L^2}^2$ and 
$\|\p^2_xu(\tau)\|_{L^2}
\leq M_0(T)(1+\tau)^{-{n\over8}-{{1}\over2}}$ 
for $s-1-\sigma(2,n)\geq 0$, namely, $s\geq[{{n-1}\over4}]+4$. 
Therefore we have 
\begin{equation}\label{I-11}
\begin{split}
I_{11}
&\leq C\!\int^{{t\over2}}_0(1+t-\tau)^{-{n\over8}-{{k+h+2}\over4}}
  \|\p^2_xu(\tau)\|_{L^2}^2d\tau \\
&\leq CM_0(T)^2\!\int^{{t\over2}}_0
  (1+t-\tau)^{-{n\over8}-{{k+h+2}\over4}}
  (1+\tau)^{-{n\over4}-1}d\tau \\
&\leq CM_0(T)^2(1+t)^{-\frac{n}{8}-\frac{k}{4}-\frac{1}{2}},
\end{split}
\end{equation}
provided that $s\geq[{{n-1}\over4}]+4$. 
On the other hand, we see that 
$\|\p^{k+h+l}_xg(\p^2_xu)\|_{L^2}\leq 
C\|\p^2_xu\|_{L^{\infty}}\|\p^{k+h+l+2}_xu\|_{L^2}$ by 
Lemma \ref{31}. Moreover, we have 
$\|\p^2_xu(\tau)\|_{L^{\infty}}\leq N_d(T)(1+\tau)^{-d}$ 
and $\|\p^{k+h+l+2}_xu(\tau)\|_{L^2}\leq E(T)$ if 
$k+h+l+2\leq s+1$, {\it i.e.}, $k+h+l\leq s-1$. Thus we obtain 
\begin{equation}\label{I-12}
\begin{split}
I_{12}
&\leq C\!\int^{{t\over2}}_0(1+t-\tau)^{-{{l+1}\over2}}
  \|\p^2_xu(\tau)\|_{L^{\infty}}
  \|\p^{k+h+l+2}_xu(\tau)\|_{L^2}d\tau \\
&\leq CN_d(T)E(T)\int^{{t\over2}}_0
  (1+t-\tau)^{-{{l+1}\over2}}(1+\tau)^{-d}d\tau
\end{split}
\end{equation} 
for $k+h+l\leq s-1$, where $d>d(n)={n\over8}+{1\over2}$. 
To obtain the optimal decay rate, we choose $l$ as the smallest 
integer satisfying 
\begin{eqnarray}\label{choice-1}
{l+1\over2}\geq\left\{
\bal
{n\over8}+{k\over4}+1-d(n)\ 
&{\rm if}\ d(n)<1,\\[1mm]
{n\over8}+{k\over4}\ &{\rm if}\ d(n)\geq1.
\ea\right.
\end{eqnarray}
This gives the choice $l=\sigma(k,n)-k$. 
In fact, when $n=2,3$, we have 
$d(n)={n\over8}+{1\over2}<1$, so that 
${l+1\over2}\geq{k\over4}+{1\over2}$. This implies 
$l\geq{k\over2}$ and hence $l=[{k+1\over2}]=\sigma_0(k)-k$. 
On the other hand, when $n\geq 4$, we have $d(n)\geq 1$, so that 
${l+1\over2}\geq{n\over8}+{k\over4}$. This implies 
$l\geq{n+2k\over4}-1$ and hence 
$l=[{n+2k-1\over4}]=\sigma_1(k,n)-k$. Thus we have 
$l=\sigma(k,n)-k$. For this choice of $l$, we obtain 
\begin{equation}\label{I-12a}
\begin{split}
I_{12}
&\leq CN_d(T)E(T)\!\int^{{t\over2}}_0
  (1+t-\tau)^{-{{l+1}\over2}}(1+\tau)^{-d}d\tau \\
&\leq CN_d(T)E(T)(1+t)^{-{n\over8}-{k\over4}},
\end{split}
\end{equation}
provided that $h$ satisfies $0\leq h\leq s-1-\sigma(k,n)$. 

Next we consider the term $I_2$. Applying \eqref{d-1} with 
$k=h+2$, $\phi=\p^k_xg(\p^2_xu)$ and $l=0$, we have 
\begin{equation}\label{I-2}
\begin{split}
I_2
&\leq C\!\int^t_{{t\over2}}(1+t-\tau)^{-{n\over8}-{{h+2}\over4}}
  \|\p^k_xg(\p^2_xu)(\tau)\|_{L^1}d\tau \\
&+C\!\int^t_{{t\over2}}(1+t-\tau)^{-{1\over2}}
  \|\p^{k+h}_xg(\p^2_xu)(\tau)\|_{L^2}d\tau
  =:I_{21}+I_{22}.
\end{split}
\end{equation}
Here we see that $\|\p^k_xg(\p^2_xu)\|_{L^1}
\leq C\|\p^2_xu\|_{L^2}\|\p^{k+2}_xu\|_{L^2}$ by Lemma \ref{31}. 
Also we know that $\|\p^2_xu(\tau)\|_{L^2}\leq M_0(T)
(1+\tau)^{-{n\over8}-{1\over2}}$ for $s\geq[{n-1\over4}]+4$. 
Moreover, we see that 
\begin{equation}\label{4x}
\|\p^{k+2}_xu(\tau)\|_{L^2}
\leq E(T)(1+\tau)^{-{k\over4}},
\end{equation}
provided that $s\geq\sigma_0(k)+1$. In fact, when $k$ is even, 
we have $\|\p^{k+2}_xu(\tau)\|_{L^2}\leq\|\p^k_xu(\tau)\|_{H^2}
\leq E(T)(1+\tau)^{-{k\over4}}$ if $s-{3\over2}k+1\geq 2$, 
{\it i.e.}, $s\geq{3\over2}k+1=\sigma_0(k)+1$. 
Also, when $k$ is odd, we have 
$\|\p^{k+2}_xu(\tau)\|_{L^2}\leq\|\p^{k+1}_xu(\tau)\|_{H^1}
\leq E(T)(1+\tau)^{-{k+1\over4}}$ if $s-{3\over2}(k+1)+1\geq 1$, 
{\it i.e.}, $s\geq{3\over2}(k+1)=\sigma_0(k)+1$. 
Thus we have verified \eqref{4x}. 
Consequently, we obtain 
\begin{equation}\label{I-21}
\begin{split}
I_{21}
&\leq C\!\int^t_{{t\over2}}(1+t-\tau)^{-{n\over8}-{{h+2}\over4}}
  \|\p^2_xu(\tau)\|_{L^2}\|\p^{k+2}_xu(\tau)\|_{L^2}d\tau \\
&\leq CM_0(T)E(T)\int^t_{{t\over2}}
  (1+t-\tau)^{-{n\over8}-{{h+2}\over4}}
  (1+\tau)^{-{n\over8}-{{k}\over4}-{1\over2}}d\tau \\
&\leq CM_0(T)E(T)(1+t)^{-{n\over8}-{k\over4}},
\end{split}
\end{equation}
provided that $s\geq[{n-1\over4}]+4$ and $s\geq\sigma_0(k)+1$. 
On the other hand, we see that 
$\|\p^{k+h}_xg(\p^2_xu)\|_{L^2}\leq C\|\p^2_xu\|_{L^\infty}
\|\p^{k+h+2}_xu\|_{L^2}$ by Lemma \ref{31}. Also we have 
$\|\p^2_xu(\tau)\|_{L^\infty}\leq N_d(T)(1+\tau)^{-d}$. 
Moreover, as a simple modification of \eqref{4x}, we have 
\begin{equation}\label{4y}
\|\p^{k+h+2}_xu(\tau)\|_{L^2}
\leq E(T)(1+\tau)^{-{k\over4}},
\end{equation}
provided that $s\geq\sigma_0(k)+h+1$. Therefore we obtain 
\begin{equation}\label{I-22}
\begin{split}
I_{22}
&\leq C\int^t_{{t\over2}}(1+t-\tau)^{-{{1}\over2}}
  \|\p^2_xu(\tau)\|_{L^{\infty}}\|\p^{k+h+2}_xu(\tau)\|_{L^2}d\tau \\
&\leq CN_d(T)E(T)\int^t_{{t\over2}}
  (1+t-\tau)^{-{{1}\over2}}(1+\tau)^{-d-{{k}\over4}}d\tau \\
&\leq CN_d(T)E(T)(1+t)^{-{n\over8}-{k\over4}}
\end{split}
\end{equation} 
for $h$ satisfying $0\leq h\leq s-1-\sigma_0(k)$, where we have 
used the requirement $d>d(n)={n\over8}+{1\over2}$. 

We substitute all these estimates into \eqref{N-0} and add 
the result for $h$ with $0\leq h\leq s-1-\sigma(k,n)$. 
This yields 
\begin{equation*}
\begin{split}
&\|\p^k_xF(u)(t)\|_{H^{s-1-\sigma(k,n)}} \\[1mm]
&\leq C\{M_0(T)^2+(M_0(T)+N_d(T))E(T)\}(1+t)^{-{n\over8}-{k\over4}}
\end{split}
\end{equation*}
for $k$ with $\sigma(k,n)\leq s-1$, where we have assumed that 
$s\geq[{n-1\over4}]+4$. This estimate together with \eqref{dL1} 
gives 
\begin{equation*}
\begin{split}
&(1+t)^{{n\over8}+{k\over4}}
  \|\p^{k}_xu(t)\|_{H^{s-1-\sigma(k,n)}} \\[1mm]
&\leq CE_1+CM_0(T)^2+C(M_0(T)+N_d(T))E(T)
\end{split}
\end{equation*}
for $k$ with $\sigma(k,n)\leq s-1$. Thus we have proved the desired 
estimate \eqref{dM-0}. 

\medskip


Second, we prove \eqref{dM-1}. We apply $\p_t$ to $F(u)$ in 
\eqref{4.3} to get
\begin{equation*}
\p_tF(u)(t)=\int^t_0\p_tG(t-\tau)\ast(1-\Delta)^{-1}
  \p^2_xg(\p^2_xu)(\tau)d\tau.
\end{equation*}
Moreover, applying $\p^{k+h}_x$ and taking the $L^2$ norm, 
we obtain 
\begin{equation}\label{N-1}
\begin{split}
&\|\p^{k+h}_x\p_tF(u)(t)\|_{L^2} \\
&\leq\int_0^t\|\p^{k+h+2}_x\p_tG(t-\tau)\ast(1-\Delta)^{-1}
  g(\p^2_xu)(\tau)\|_{L^2}d\tau \\
&=\int^{t\over2}_0+\int^t_{t\over2}=:I_3+I_4.
\end{split}
\end{equation}
We estimate the term $I_3$ by applying \eqref{d-2} with $k$ 
replaced by $k+h+2$ and with $\phi=g(\p^2_xu)$ as 
\begin{equation}\label{I-3}
\begin{split}
I_3
&\leq C\!\int^{{t\over2}}_0(1+t-\tau)^{-{n\over8}-{{k+h+2}\over4}-1}
  \|g(\p^2_xu)(\tau)\|_{L^1}d\tau \\
&+C\!\int^{{t\over2}}_0(1+t-\tau)^{-{{l}\over2}}
  \|\p^{k+h+l}_xg(\p^2_xu)(\tau)\|_{L^2}d\tau
  =:I_{31}+I_{32}.
\end{split}
\end{equation}
Here the term $I_{31}$ is estimated just in the same way as 
$I_{11}$ in \eqref{I-11} and we have 
\begin{equation}\label{I-31}
\begin{split}
I_{31}
&\leq C\!\int^{{t\over2}}_0(1+t-\tau)^{-{n\over8}-{{k+h+2}\over4}-1}
  \|\p^2_xu(\tau)\|_{L^2}^2d\tau \\
&\leq CM_0(T)^2\!\int^{{t\over2}}_0
  (1+t-\tau)^{-{n\over8}-{{k+h+2}\over4}-1}
  (1+\tau)^{-{n\over4}-1}d\tau \\
&\leq CM_0(T)^2(1+t)^{-\frac{n}{8}-\frac{k}{4}-\frac{3}{2}},
\end{split}
\end{equation}
provided that $s\geq[{{n-1}\over4}]+4$. 
Also, similarly to $I_{12}$ in \eqref{I-12}, we get 
\begin{equation}\label{I-32}
\begin{split}
I_{32}
&\leq C\!\int^{{t\over2}}_0(1+t-\tau)^{-{{l}\over2}}
  \|\p^2_xu(\tau)\|_{L^{\infty}}
  \|\p^{k+h+l+2}_xu(\tau)\|_{L^2}d\tau \\
&\leq CN_d(T)E(T)\int^{{t\over2}}_0
  (1+t-\tau)^{-{{l}\over2}}(1+\tau)^{-d}d\tau
\end{split}
\end{equation}
for $k+h+l\leq s-1$, where $d>d(n)={n\over8}+{1\over2}$. 
We choose $l$ as the smallest integer satisfying 
\begin{eqnarray}\label{choice-2}
{l\over2}\geq\left\{
\bal
{n\over8}+{k\over4}+2-d(n)\ 
&{\rm if}\ d(n)<1,\\[1mm]
{n\over8}+{k\over4}+1\ &{\rm if}\ d(n)\geq1.
\ea\right.
\end{eqnarray}
A similar observation as in \eqref{choice-1} shows that the desired 
choice is $l=\sigma(k,n)-k+3$. For this choice of $l$, we obtain 
\begin{equation}\label{I-32a}
\begin{split}
I_{32}
&\leq CN_d(T)E(T)\!\int^{{t\over2}}_0
  (1+t-\tau)^{-{{l}\over2}}(1+\tau)^{-d}d\tau \\
&\leq CN_d(T)E(T)(1+t)^{-{n\over8}-{k\over4}-1}
\end{split}
\end{equation}
for $h$ satisfying $0\leq h\leq s-4-\sigma(k,n)$; here we need to 
assume that $s\geq \sigma(0,n)+4=[{n-1\over4}]+4$. 

For the term $I_4$, we apply \eqref{d-2} with $k=h$, 
$\phi=\p^{k+2}_xg(\p^2_xu)$ and $l=2$. This gives 
\begin{equation}\label{I-4}
\begin{split}
I_4
&\leq C\!\int^t_{{t\over2}}(1+t-\tau)^{-{n\over8}-{{h}\over4}-1}
  \|\p^{k+2}_xg(\p^2_xu)(\tau)\|_{L^1}d\tau \\
&+C\!\int^t_{{t\over2}}(1+t-\tau)^{-1}
  \|\p^{k+h+2}_xg(\p^2_xu)(\tau)\|_{L^2}d\tau
  =:I_{41}+I_{42}.
\end{split}
\end{equation}
Here we have $\|\p^{k+2}_xg(\p^2_xu)\|_{L^1}
\leq C\|\p^2_xu\|_{L^2}\|\p^{k+4}_xu\|_{L^2}$ by Lemma \ref{31}. 
Also we know that $\|\p^2_xu(\tau)\|_{L^2}\leq M_0(T)
(1+\tau)^{-{n\over8}-{1\over2}}$ for $s\geq[{n-1\over4}]+4$. 
Moreover, as a counterpart of \eqref{4x}, we get 
$\|\p^{k+4}_xu(\tau)\|_{L^2}\leq E(T)(1+\tau)^{-{k+2\over4}}$,
provided that $s\geq\sigma_0(k+2)+1=\sigma_0(k)+4$. 
Therefore we obtain 
\begin{equation}\label{I-41}
\begin{split}
I_{41}
&\leq C\!\int^t_{{t\over2}}(1+t-\tau)^{-{n\over8}-{{h}\over4}-1}
  \|\p^2_xu(\tau)\|_{L^2}\|\p^{k+4}_xu(\tau)\|_{L^2}d\tau \\
&\leq CM_0(T)E(T)\int^t_{{t\over2}}
  (1+t-\tau)^{-{n\over8}-{{h}\over4}-1}
  (1+\tau)^{-{n\over8}-{{k+2}\over4}-{1\over2}}d\tau \\
&\leq CM_0(T)E(T)(1+t)^{-{n\over8}-{k\over4}-1},
\end{split}
\end{equation}
provided that $s\geq[{n-1\over4}]+4$ and $s\geq\sigma_0(k)+4$. 
On the other hand, we see that 
$\|\p^{k+h+2}_xg(\p^2_xu)\|_{L^2}\leq C\|\p^2_xu\|_{L^\infty}
\|\p^{k+h+4}_xu\|_{L^2}$ by Lemma \ref{31}. Also, as a conterpart 
of \eqref{4y}, we get 
$\|\p^{k+h+4}_xu(\tau)\|_{L^2}\leq E(T)(1+\tau)^{-{k+2\over4}}$, 
provided that $s\geq\sigma_0(k+2)+h+1=\sigma_0(k)+h+4$. 
Therefore, similarly to $I_{22}$ in \eqref{I-22}, we obtain 
\begin{equation}\label{I-42}
\begin{split}
I_{42}
&\leq C\int^t_{{t\over2}}(1+t-\tau)^{-1}
  \|\p^2_xu(\tau)\|_{L^{\infty}}\|\p^{k+h+4}_xu(\tau)\|_{L^2}d\tau \\
&\leq CN_d(T)E(T)\int^t_{{t\over2}}
  (1+t-\tau)^{-1}(1+\tau)^{-d-{{k+2}\over4}}d\tau \\
&\leq CN_d(T)E(T)(1+t)^{-{n\over8}-{k\over4}-1}
\end{split}
\end{equation}
for $h$ satisfying $0\leq h\leq s-4-\sigma_0(k)$, where we have 
used the requirement $d>d(n)={n\over8}+{1\over2}$. 

Substituting all these estimates into \eqref{N-1} and adding 
for $h$ with $0\leq h\leq s-4-\sigma(k,n)$, we arrive at 
\begin{equation*}
\begin{split}
&\|\p^k_x\p_tF(u)(t)\|_{H^{s-4-\sigma(k,n)}} \\[1mm]
&\leq C\{M_0(T)^2+(M(T)+N_d(T))E(T)\}(1+t)^{-{n\over8}-{k\over4}-1}
\end{split}
\end{equation*}
for $k$ with $\sigma(k,n)\leq s-4$, where we have assumed that 
$s\geq[{n-1\over4}]+4$. This estimate together with \eqref{dL2} 
gives the desired estimate \eqref{dM-1}. 
Therefore the proof of Proposition \ref{p42} is complete. 
\end{proof}


\section{Proof of Theorem \ref{21}}

The aim of this section is to prove Theorem \ref{21}. 
Since a local existence result is obtained in Theorem \ref{LS}, 
we only need to show the a priori estimates stated in Proposition 
\ref{AP}. 

To complete the estimates for the time-weighted $L^2$ norms 
$E(T)$, $M_0(T)$ and $M_1(T)$, we need to control the special 
$L^\infty$ norms $L(T)$ and $N_d(T)$ with 
$d>d(n)={n\over8}+{1\over2}$. This is done by the following lemma. 

\bl\label{42}
Let $n\geq 2$. We have 
\begin{equation}\label{L}
L(T)\leq C(M_0(T)+M_1(T)+E(T))
\end{equation}
if $s\geq 8$ for $n=2$ and $s\geq[{n\over2}]+5$ for $n\geq 3$. 
Also we have 
\begin{equation}\label{Nd}
N_d(T)\leq C(M_0(T)+E(T))
\end{equation}
for some $d$ with $d>d(n)={n\over8}+{1\over2}$, provided that 
$s\geq 3[{n\over4}]+5$ for $n\geq 2$. 
\el

\begin{proof}
For the proof of \eqref{L}, it suffices to show the following 
decay estimates: 
\begin{eqnarray}
&\|\p^3_xu(t)\|_{L^\infty}\leq C((M_0(T)+E(T))(1+t)^{-\gamma}, 
\label{dL-1}
\\[1mm]
&\|\p^2_xu_t(t)\|_{L^\infty}
  \leq C((M_0(T)+M_1(T)+E(T))(1+t)^{-\gamma}
\label{dL-2}
\end{eqnarray}
for some $\gamma>1$. These decay estimates are proved as follows. 
First, applying the Gagliado-Nirenberg inequality, we have 
\begin{equation}\label{5.5}
\begin{split}
&\|\p^3_xu\|_{L^\infty}
  \leq C\|\p^3_xu\|_{L^2}^{1-\theta}
  \|\p^{s_0+3}_xu\|_{L^2}^\theta, \\[1mm]
&\|\p^2_xu_t\|_{L^\infty}
  \leq C\|\p^2_xu_t\|_{L^2}^{1-\theta}
  \|\p^{s_0+2}_xu_t\|_{L^2}^\theta,
\end{split}
\end{equation}
where $s_0=[{n\over2}]+1$ and $\theta={n\over2s_0}<1$. 
Here we see that 
\begin{equation}\label{5.6}
\|\p^3_xu(t)\|_{L^2}\leq M_0(T)(1+t)^{-{n\over8}-{3\over4}},
\end{equation}
provided that $s-1-\sigma(3,n)\geq 0$, that is, $s\geq 6$ for 
$n=2,\,3$ and $s\geq[{n+1\over4}]+5$ for $n\geq 3$. 
Notice that the decay rate here verifies ${n\over8}+{3\over4}>1$ for 
$n\geq3$ and ${n\over8}+{3\over4}=1$ for $n=2$. We see that 
$$
\|\p^{s_0+3}_xu(t)\|_{L^2}\leq \|\p^{4}_xu(t)\|_{H^{s_0-1}}
\leq E(T)(1+t)^{-1}
$$
if $s-3\times2+1\geq s_0-1$, {\it i.e.}, $s\geq[{n\over2}]+5$. 
Substituing this estimate and \eqref{5.6} into the first inequality 
in \eqref{5.5}, we obtain \eqref{dL-1} with 
$\gamma=({n\over8}+{3\over4})(1-\theta)+\theta$ for 
$s\geq[{n\over2}]+5$. Since this decay rate $\gamma$ verifies 
$\gamma=1+({n\over8}-{1\over4})(1-\theta)>1$ for $n\geq3$, 
we have shown \eqref{dL-1} for $s\geq[{n\over2}]+5$ and $n\geq3$. 

Now we restrict our attention to the case $n=2$ in which 
${n\over8}+{3\over4}=1$. 
In this case the first inequality in \eqref{5.5} becomes 
$$
\|\p^3_xu\|_{L^\infty}\leq C\|\p^3_xu\|_{L^2}^{1\over2}
\|\p^{5}_xu\|_{L^2}^{1\over2}. 
$$
Here, using a simple interpolation inequality 
$\|\p_xv\|_{L^2}\leq C\|v\|_{L^2}^{1\over2}
\|\p^2_xv\|_{L^2}^{1\over2}$ with $v=\p^{4}_xu$, we have 
$$
\|\p^5_xu(t)\|_{L^2}\leq CE(T)(1+t)^{-{5\over4}}
$$
for $s-3\times3+1\geq0$, {\it i.e.}, $s\geq8$. Thus, even for $n=2$, 
we get the desired estimate \eqref{dL-1} with $\gamma={9\over8}$ 
for $s\geq8$. 

The proof of \eqref{dL-2} is similar. We see that 
$$
\|\p^2_xu_t(t)\|_{L^2}\leq \|u_t(t)\|_{H^2} 
\leq M_1(T)(1+t)^{-{n\over8}-1}
$$
if $s-4-\sigma(0,n)\geq2$, {\it i.e.}, $s\geq[{n-1\over4}]+6$. 
Also, we have 
$$
\|\p^{s_0+2}_xu_t(t)\|_{L^2}\leq \|\p^{2}_xu_t(t)\|_{H^{s_0}}
\leq E(T)(1+t)^{-1},
$$
provided that $s-3\times1-1\geq s_0$, {\it i.e.}, 
$s\geq[{n\over2}]+5$. These estimates together with the second 
inequality in \eqref{5.5} prove \eqref{dL-2} with 
$\gamma=({n\over8}+1)(1-\theta)+\theta=1+{n\over8}(1-\theta)>1$ 
for $s\geq[{n\over2}]+5$. 

Finally, we prove \eqref{Nd}. 
Applying the Gagliardo-Nirenberg inequlity, we have 
\begin{equation}\label{5.7}
\|\p^2_xu\|_{L^\infty}
  \leq C\|\p^2_xu\|_{L^2}^{1-\theta}
  \|\p^{s_0+2}_xu\|_{L^2}^\theta,
\end{equation}
where $s_0$ and $\theta$ are the same as in \eqref{5.5}. 
Here we see that 
$$
\|\p^2_xu(t)\|_{L^2}\leq M_0(T)(1+t)^{-{n\over8}-{1\over2}}
$$
if $s-1-\sigma(2,n)\geq 0$, {\it i.e.}, $s\geq[{n-1\over4}]+4$. 
Moreover, we can show that 
\begin{equation}\label{5x}
\|\p^{s_0+2}_xu(t)\|_{L^2}\leq CE(T)(1+t)^{-{s_0+2\over4}},
\end{equation}
provided that $s\geq 3[{n\over4}]+5$. 
In fact, when $s_0=[{n\over2}]+1$ is odd ({\it i.e.}, $[{n\over2}]$ 
is even), using an interpolation inequality 
$\|\p_xv\|_{L^2}\leq C\|v\|_{L^2}^{1\over2}
\|\p^2_xu\|_{L^2}^{1\over2}$ with $v=\p^{s_0+1}_xu$ and noting the 
definition of $E(T)$, we obtain the estimate \eqref{5x}, 
provided that $s-{3\over2}(s_0+3)+1\geq 0$, {\it i.e.}, 
$s\geq{3\over2}[{n\over2}]+5$. Since $[{n\over2}]$ is even 
in this case, this requirement on $s$ can be rewritten as 
$s\geq 3[{n\over4}]+5$. 
On the other hand, when $s_0=[{n\over2}]+1$ is even 
({\it i.e.}, $[{n\over2}]$ is odd), we easily get \eqref{5x}, 
provided that $s-{3\over2}(s_0+2)+1\geq 0$, {\it i.e.}, 
$s\geq{3\over2}[{n\over2}]+{7\over2}$. 
Since $[{n\over2}]$ is odd in this case, 
this requirement on $s$ is also rewritten as $s\geq 3[{n\over4}]+5$. 
Thus we have shown the estimate \eqref{5x}. 
Substituting these estimates into \eqref{5.7}, we obatin 
$$
\|\p^2_xu(t)\|\leq C(M_0(T)+E(T))(1+t)^{-d}
$$
with $d=({n\over8}+{1\over2})(1-\theta)+{s_0+2\over4}\,\theta$, 
provided that $s\geq 3[{n\over4}]+5$ and $n\geq2$. 
Note that this decay rate $d$ verifies 
$d=({n\over8}+{1\over2})+{1\over4}(s_0-{n\over2})\theta
>{n\over8}+{1\over2}$ because $s_0=[{n\over2}]+1>{n\over2}$. 
This proves the desired estimate \eqref{Nd}. 
Thus the proof of Lemma \ref{42} is complete. 
\end{proof}

\begin{proof}[Proof of Proposition \ref{AP} and Theorem \ref{21}]\ 
Put 
$$
Y(T):=E(T)+D(T)+M_0(T)+M_1(T).
$$
By virtue of Propositions \ref{p41} and \ref{p42} together with 
Lemma \ref{42}, we have the inequlity $Y(T)^2\leq CE_1^2+CY(T)^3$, 
from which we can deduce that $Y(T)\leq CE_1$, provided that 
$E_1$ is suitably small, namely, $E_1\leq\d_1$, where $\d_1$ is a 
small positive number depending on $\bar{\d}$ in \eqref{3a}. 
This gives the desired a priori estimates \eqref{energy} and 
\eqref{decay} for solutions satisfying \eqref{3a}. Thus we have 
proved Proposition \ref{AP}. 

By virtue of the a priori estimate \eqref{energy}, we can continue 
a unique local solution obtained in Theorem \ref{LS} 
globally in time, provided that $E_1$ is suitably small, say, 
$E_1\leq\delta_0$. The global solution thus obtained satisfies 
\eqref{energy} and \eqref{decay} for any $T>0$. 
In particular, we have the decay estimates \eqref{decay-1} and 
\eqref{decay-2} from \eqref{decay}. 
This completes the proof of Theorem \ref{21}. 
\end{proof}


\section{Asymptotic profile}

The aim of this section is to prove Theorem \ref{22} on the 
asymptotic profile. First we prove that the solution to the problem 
\eqref{1a}, \eqref{IC} can be approximated by the solution to the 
corresponding linear problem \eqref{1cc}, \eqref{IC}. 

\bp\label{p51} 
Suppose that {\rm [A1]} and {\rm [A2]} are 
satisfied. Let $n\geq 2$ and $s\geq s(n)$. Assume that 
$u_0\in H^{s+1}(\R^n)\cap L^1(\R^n)$ and 
$u_1\in H^s(\R^n)\cap L^1(\R^n)$, 
and put $E_1=\|u_0\|_{H^{s+1}}+\|u_1\|_{H^s}+\|(u_0,u_1)\|_{L^1}$. 
Let $u(x,t)$ be the global solution to the problem \eqref{1a}, 
\eqref{IC} which is obtained in Theorem \ref{21}, and let 
$\bar{u}(x,t)$ be the solution to the corresponding linear problem 
which is given by the formula \eqref{4.3}. Then we have 
\begin{eqnarray}
&\|\p^k_x(u-\bar{u})(t)\|_{H^{s-2-\sigma(k,n)}}
  \leq CE_1^2(1+t)^{-{n\over8}-{k\over4}-{1\over2}},
\label{6.1}
\\[1mm]
&\|\p^k_x\p_t(u-\bar{u})(t)\|_{H^{s-6-\sigma(k,n)}}
  \leq CE_1^2(1+t)^{-{n\over8}-{k\over4}-{3\over2}}
\label{6.2}
\end{eqnarray}
for $k\geq 0$, where $\sigma(k,n)\leq s-2$ in \eqref{6.1} and 
$\sigma(k,n)\leq s-6$ in \eqref{6.2}. 
\ep

\begin{proof}
It suffices to show the following estimates for the nonlinear part 
$F(u)(x,t)$ given in \eqref{4.3}. 
\begin{eqnarray}
&\|\p^k_xF(u)(t)\|_{H^{s-2-\sigma(k,n)}}
  \leq CE_1^2(1+t)^{-{n\over8}-{k\over4}-{1\over2}},
\label{6.3}
\\[1mm]
&\|\p^k_x\p_tF(u)(t)\|_{H^{s-6-\sigma(k,n)}}
  \leq CE_1^2(1+t)^{-{n\over8}-{k\over4}-{3\over2}}.
\label{6.4}
\end{eqnarray}
We first prove \eqref{6.3} similarly as in the proof of \eqref{dM-0}. 
We have \eqref{N-0} and \eqref{I-1}. Moreover, the term $I_{11}$ in 
\eqref{I-1} is estimated just in the same way as in \eqref{I-11} 
and we obtain 
\begin{equation*}
I_{11} \leq CE_1^2(1+t)^{-\frac{n}{8}-\frac{k}{4}-\frac{1}{2}}
\end{equation*}
for $s\geq[{n-1\over4}]+4$. 
For the term $I_{12}$ in \eqref{I-1}, we also have \eqref{I-12} for 
$k+h+l\leq s-1$. 
To get a better decay estimate, we choose $l$ in \eqref{I-12} as the 
smallest integer satisfying 
\begin{eqnarray*} 
{l+1\over2}\geq\left\{ \bal {n\over8}+{k\over4}+{3\over2}-d(n)\
&{\rm if}\ d(n)<1,\\[1mm]
{n\over8}+{k\over4}+{1\over2}\ &{\rm if}\ d(n)\geq1. \ea\right.
\end{eqnarray*}
This gives the choice $l=\sigma(k,n)-k+1$. For this choice of $l$, 
as a counterpart of \eqref{I-12a}, we obtain 
\begin{equation*} 
I_{12} \leq CE_1^2\!\int^{{t\over2}}_0
  (1+t-\tau)^{-{{l+1}\over2}}(1+\tau)^{-d}d\tau 
\leq CE_1^2(1+t)^{-{n\over8}-{k\over4}-{1\over2}}
\end{equation*}
for $h$ with $0\leq h\leq s-2-\sigma(k,n)$. 

Next we estimate the term $I_2$ in \eqref{N-0}. 
By applying \eqref{d-1} with $k=h+2$ and $\phi=\p^k_xg(\p^2_xu)$, 
as a counterpart of \eqref{I-2}, we have 
\begin{equation}\label{6I-2}
\begin{split}
I_2 &\leq C\!\int^t_{{t\over2}}(1+t-\tau)^{-{n\over8}-{{h+2}\over4}}
  \|\p^k_xg(\p^2_xu)(\tau)\|_{L^1}d\tau \\
&+C\!\int^t_{{t\over2}}(1+t-\tau)^{-{l+1\over2}}
  \|\p^{k+h+l}_xg(\p^2_xu)(\tau)\|_{L^2}d\tau
  =:I_{21}+I^{\prime}_{22},
\end{split}
\end{equation}
where $l+1\geq 0$. 
To show a better decay estimate for $I_{21}$, instead of \eqref{4x}, 
we use the following estimates: 
\begin{eqnarray}
&\|\p^{k+2}_xu(\tau)\|_{L^2} \leq CE_1(1+\tau)^{-{k+2\over4}},
\label{64x}
\\[1mm]
&\|\p^{k+2}_xu(\tau)\|_{L^2} \leq CE_1(1+\tau)^{-{k+1\over4}},
\label{6-4x}
\end{eqnarray}
where $k$ is even and $s\geq\sigma_0(k)+2$ in \eqref{64x}, and 
$k$ is odd and $s\geq\sigma_0(k)+1$ in \eqref{6-4x}. 
These estimates are verified as follows.
When $k$ is even, we have 
$\|\p^{k+2}_xu(\tau)\|_{L^2} \leq C E_1(1+\tau)^{-{k+2\over4}}$ 
if $s-{3\over2}(k+2)+1\geq 0$, {\it i.e.}, 
$s\geq{3\over2}k+2=\sigma_0(k)+2$. Also, when $k$ is odd, we have 
$\|\p^{k+2}_xu(\tau)\|_{L^2}\leq\|\p^{k+1}_xu(\tau)\|_{H^1} \leq
CE_1(1+\tau)^{-{k+1\over4}}$ if $s-{3\over2}(k+1)+1\geq 1$, 
{\it i.e.}, $s\geq{3\over2}(k+1)=\sigma_0(k)+1$. 
Thus we have verified \eqref{64x} and \eqref{6-4x}. 
Now, using \eqref{6-4x} for $s\geq\sigma_0(k)+2$, we obtain 
\begin{equation*} 
\begin{split}
I_{21} &\leq
C\!\int^t_{{t\over2}}(1+t-\tau)^{-{n\over8}-{{h+2}\over4}}
  \|\p^2_xu(\tau)\|_{L^2}\|\p^{k+2}_xu(\tau)\|_{L^2}d\tau \\
&\leq CE_1^2\!\int^t_{{t\over2}}
  (1+t-\tau)^{-{n\over8}-{{h+2}\over4}}
  (1+\tau)^{-{n\over8}-{{k+1}\over4}-{1\over2}}d\tau \\
&\leq CE_1^2(1+t)^{-{n\over8}-{k\over4}-{1\over2}},
\end{split}
\end{equation*}
provided that $s\geq[{n-1\over4}]+4$ and $s\geq\sigma_0(k)+2$, 
where we have used $n\geq 2$. 

Finally, we estimate $I^{\prime}_{22}$. 
We choose $l$ in $I^{\prime}_{22}$ such that $l=0$ if $k$ is even 
and $l=1$ if $k$ is odd. Also, we see that 
\begin{eqnarray}
&\|\p^{k+h+2}_xu(\tau)\|_{L^2} \leq CE_1(1+\tau)^{-{k+2\over4}},
\label{64y}
\\[1mm]
&\|\p^{k+h+3}_xu(\tau)\|_{L^2} \leq CE_1(1+\tau)^{-{k+1\over4}}
\label{6-4y}
\end{eqnarray}
for $s\geq\sigma_0(k)+h+2$, where $k$ is even in \eqref{64y} and
$k$ is odd in \eqref{6-4y}; these are simple modifications of 
\eqref{64x} and \eqref{6-4x}, respectively. 
Therefore, when $k$ is even, by taking $l=0$ in $I^{\prime}_{22}$ 
and using \eqref{64y}, we obtain 
\begin{equation*} 
\begin{split}
I^{\prime}_{22} &\leq C\!\int^t_{{t\over2}}(1+t-\tau)^{-{1\over2}}
  \|\p^2_xu(\tau)\|_{L^{\infty}}
  \|\p^{k+h+2}_xu(\tau)\|_{L^2}d\tau \\
&\leq CE_1^2\!\int^t_{{t\over2}}(1+t-\tau)^{-{1\over2}}
  (1+\tau)^{-d-{k+2\over4}}d\tau 
  \leq CE_1^2(1+t)^{-{n\over8}-{k\over4}-{1\over2}}
\end{split}
\end{equation*}
for $h$ with $0\leq h\leq s-2-\sigma_0(k)$, where we have used 
the requirement $d>d(n)={n\over8}+{1\over2}$. 
On the other hand, when $k$ is odd, by taking $l=1$ and using 
\eqref{6-4y}, we get 
\begin{equation*} 
\begin{split}
I^{\prime}_{22} &\leq C\!\int^t_{{t\over2}}(1+t-\tau)^{-1}
  \|\p^2_xu(\tau)\|_{L^{\infty}}
  \|\p^{k+h+3}_xu(\tau)\|_{L^2}d\tau \\
&\leq CE_1^2\!\int^t_{{t\over2}}(1+t-\tau)^{-1}
  (1+\tau)^{-d-{k+1\over4}}d\tau 
  \leq CE_1^2(1+t)^{-{n\over8}-{k\over4}-{1\over2}}
\end{split}
\end{equation*}
for $h$ with $0\leq h\leq s-2-\sigma_0(k)$, where we have again 
used $d>d(n)={n\over8}+{1\over2}$. 

We substitute all these estimates into \eqref{N-0} and add the 
result for $h$ with $0\leq h\leq s-2-\sigma(k,n)$. This yields 
the desired estimate \eqref{6.3}. 

\medskip

Second, we prove \eqref{6.4} similarly as in the proof of 
\eqref{dM-1}. In this case we also have \eqref{N-1} and \eqref{I-3}. 
Here the term $I_{31}$ in \eqref{I-3} is estimated in the same way 
as in \eqref{I-31} and we have 
\begin{equation*} 
\begin{split}
I_{31} \leq CE_1^2(1+t)^{-\frac{n}{8}-\frac{k}{4}-\frac{3}{2}}
\end{split}
\end{equation*}
for $s\geq[{n-1\over4}]+4$. 
For the term $I_{32}$ in \eqref{I-3}, we also have \eqref{I-32} for 
$k+h+l\leq s-1$. 
To get a better decay estimate, we choose $l$ as the smallest 
integer satisfying 
\begin{eqnarray*} 
{l\over2}\geq\left\{ \bal {n\over8}+{k\over4}+{5\over2}-d(n)\
&{\rm if}\ d(n)<1,\\[1mm]
{n\over8}+{k\over4}+{3\over2}\ &{\rm if}\ d(n)\geq1. \ea\right.
\end{eqnarray*}
This gives $l=\sigma(k,n)-k+4$. For this choice of $l$, 
as a counterpart of \eqref{I-32a}, we obtain 
\begin{equation*} 
I_{32} \leq CE_1^2\!\int^{{t\over2}}_0
  (1+t-\tau)^{-{{l}\over2}}(1+\tau)^{-d}d\tau 
  \leq CE_1^2(1+t)^{-{n\over8}-{k\over4}-{3\over2}}
\end{equation*}
for $h$ with $0\leq h\leq s-5-\sigma(k,n)$; here we need to 
assume that $s\geq \sigma(0,n)+5=[{n-1\over4}]+5$. 

For the term $I_4$, we have \eqref{I-4}. To estimate the term 
$I_{41}$ in \eqref{I-4}, we use the estimate 
\begin{equation}\label{6.10}
\|\p^{k+4}_xu(\tau)\|_{L^2}\leq CE_1(1+\tau)^{-{k\over4}-1}
\end{equation}
for $s\geq\sigma(k,n)+6$. This is verified as follows. 
When $k$ even, we have 
$\|\p^{k+4}_xu(\tau)\|_{L^2}\leq CE_1(1+\tau)^{-{k\over4}-1}$ 
if $s-{3\over2}(k+4)+1\geq0$, {\it i.e.}, $s\geq \sigma_0(k)+5$. 
On the other hand, when $k$ is odd and $n=2,3$, we see that 
$\|\p^{k+4}_xu(\tau)\|_{L^2}\leq\|\p^{k+3}_xu(\tau)\|_{H^1}
\leq CE_1(1+\tau)^{-{n\over8}-{k+3\over4}}$ if 
$s-1-\sigma_0(k+3)\geq 1$, {\it i.e.}, $s\geq\sigma_0(k)+6$. 
Also, when $k$ is odd and $n\geq4$, we have 
$\|\p^{k+4}_xu(\tau)\|_{L^2}\leq\|\p^{k+2}_xu(\tau)\|_{H^2}
\leq CE_1(1+\tau)^{-{n\over8}-{k+2\over4}}$ if 
$s-1-\sigma(k+2,n)\geq 2$, {\it i.e.}, $s\geq\sigma(k,n)+6$. 
These considerations prove \eqref{6.10}. 
By using \eqref{6.10}, as a counterpart of \eqref{I-41}, we obtain 
\begin{equation*} 
\begin{split}
I_{41} &\leq
C\!\int^t_{{t\over2}}(1+t-\tau)^{-{n\over8}-{{h}\over4}-1}
  \|\p^2_xu(\tau)\|_{L^2}\|\p^{k+4}_xu(\tau)\|_{L^2}d\tau \\
&\leq CE_1^2\!\int^t_{{t\over2}}
  (1+t-\tau)^{-{n\over8}-{{h}\over4}-1}
  (1+\tau)^{-{n\over8}-{{k}\over4}-{3\over2}}d\tau \\
&\leq CE_1^2(1+t)^{-{n\over8}-{k\over4}-{3\over2}},
\end{split}
\end{equation*}
provided that $s\geq\sigma(k,n)+6$. 

Finally, we estiamte the term $I_{42}$ in \eqref{I-4}. 
Similarly as in \eqref{6.10}, we have 
$\|\p^{k+h+4}_xu(\tau)\|_{L^2}\leq CE_1(1+\tau)^{-{k\over4}-1}$, 
provided that $s\geq\sigma(k,n)+h+6$. 
Therefore, as a counterpart of \eqref{I-42}, we obtain 
\begin{equation*} 
\begin{split}
I_{42} &\leq C\!\int^t_{{t\over2}}(1+t-\tau)^{-1}
  \|\p^2_xu(\tau)\|_{L^{\infty}}\|\p^{k+h+4}_xu(\tau)\|_{L^2}d\tau \\
&\leq CE_1^2\!\int^t_{{t\over2}}
  (1+t-\tau)^{-1}(1+\tau)^{-d-{{k}\over4}-1}d\tau 
\leq CE_1^2(1+t)^{-{n\over8}-{k\over4}-{3\over2}}
\end{split}
\end{equation*}
for $h$ satisfying $0\leq h\leq s-6-\sigma(k,n)$, where we have used 
the requirement $d>d(n)={n\over8}+{1\over2}$. 

Substituting all these estimates into \eqref{N-1} and adding for $h$ 
with $0\leq h\leq s-6-\sigma(k,n)$, we arrive at the desired 
estimate \eqref{6.4}. 
This completes the proof of Proposition \ref{p51}. 
\end{proof}


Next, by similar computations as in \cite{SK}, we know that 
the solution $\bar{u}(x,t)$ to the linear problem \eqref{1cc}, 
\eqref{IC} is well approximated for $t\to\infty$ by the function 
$(G_0(t)*(u_0+u_1))(x)$, which is the solution to the fourth-order 
linear parabolic equation \eqref{limit-eq} with 
the initial data $u(x,0)=u_0(x)+u_1(x)$. Here $G_0(x,t)$ is the 
fundamental solution of \eqref{limit-eq} given in \eqref{2.7}. 
More precisely, we have: 

\bl\label{Su1}
Let $n\geq1$ and assume that $u_0\in H^{s+1}(\R^n)\cap L^1(\R^n)$ 
and $u_1\in H^s(\R^n)\cap L^1(\R^n)$ for an integer $s$ specified 
below. Put 
$E_1=\|u_0\|_{H^{s+1}}+\|u_1\|_{H^s}+\|(u_0,u_1)\|_{L^1}$. 
Then the linear part $\bar{u}(x,t)$ in \eqref{4.3} satisfies the 
following asymptotic relations: 
\begin{equation}\label{6.12}
\|\p^k_x\{\bar{u}(t)-G_0(t)*(u_0+u_1)\}\|_{H^{s-1-\sigma_1(k,n)}}
  \leq CE_1(1+t)^{-{n\over8}-{k\over4}-{1\over2}}, 
\end{equation}
\begin{equation}\label{6.13}
\|\p^k_x\p_t\{\bar{u}(t)-G_0(t)*(u_0+u_1)\}\|_{H^{s-4-\sigma_1(k,n)}}
  \leq CE_1(1+t)^{-{n\over8}-{k\over4}-{3\over2}},
\end{equation}
where $s\geq[{n-1\over4}]+1$, $k\geq 0$ and $\sigma_1(k,n)\leq s-1$ 
in \eqref{6.12}, and 
$s\geq[{n-1\over4}]+4$, $k\geq 0$ and $\sigma_1(k,n)\leq s-4$ in 
\eqref{6.13}. 
Here $G_0(x,t)$ is the fundamental solution of \eqref{limit-eq} 
given in \eqref{2.7}. 
\el

The solution of \eqref{limit-eq} with the initial data 
$u(x,0)=u_0(x)+u_1(x)$ is further approximated by the self-similar 
solution $MG_0(x,t+1)=M(t+1)^{-{n\over4}}\phi_0(x/(t+1)^{{1\over4}})$, 
where $M=\int_{\mathbb{R}^n}(u_0+u_1)(x)dx$, and $\phi_0(x)$ is 
given by \eqref{2.9}. 
To see this, we need the following preparation on the solution 
operator $G_0(t)*$. 

\bl\label{51} 
Let $n\geq 1$ and $s\geq 0$. Assume that 
$\phi\in H^s(\mathbb{R}^n)\cap L^1(\mathbb{R}^n)$. 
Then we have 
\begin{equation}\label{6.14}
\begin{split}
&\|\p^k_xG_0(t)*\phi\|_{L^2}
  \leq C(1+t)^{-{n\over8}-{k\over4}}\|\phi\|_{L^1}
  +Ce^{-ct}\|\p^k_x\phi\|_{L^2}, \\[1mm]
&\|\p^k_x\p_tG_0(t)*\phi\|_{L^2}
  \leq C(1+t)^{-{n\over8}-{k\over4}-1}\|\phi\|_{L^1}
  +Ce^{-ct}\|\p^{k+4}_x\phi\|_{L^2}
\end{split}
\end{equation}
for $0\leq k\leq s$ and $0\leq k\leq s-4$, respectively; 
we have assumed $s\geq 4$ in the latter estimate. 
Also, if $\phi\in L^1_1(\mathbb{R}^n)$ and 
$\int_{\mathbb{R}^n}\phi(x)dx=0$, then we have 
\begin{equation}\label{6.15}
\begin{split}
&\|\p^k_xG_0(t)*\phi\|_{L^2}
  \leq Ct^{-{n\over8}-{{k+1}\over4}}\|\phi\|_{L^1_1}, \\[1mm]
&\|\p^k_x\p_tG_0(t)*\phi\|_{L^2}
  \leq Ct^{-{n\over8}-{{k+5}\over4}}\|\phi\|_{L^1_1}.
\end{split}
\end{equation}
\el

\begin{proof}
The proof of \eqref{6.14} is easy and we omit it. 
Here we only prove the first estimate in \eqref{6.15}; the second 
one in \eqref{6.15} is proved similarly. 
The assumption $\int_{\mathbb{R}^n}\phi(x)dx=0$ implies that 
$\hat{\phi}(0)=0$ and hence we have 
$$
|\hat{\phi}(\xi)|
\leq|\xi|\,\|\p_\xi\hat{\phi}\|_{L^\infty}
\leq|\xi|\,\|\phi\|_{L^1_1}.
$$
Therefore, applying the Plancherel theorem, we obtain 
\begin{equation*}
\begin{split}
&\|\p^k_xG_0(t)*\phi\|_{L^2}^2
  =\frac{1}{(2\pi)^n}\int_{\mathbb{R}^n}
  |\xi|^{2k}e^{-2\gamma(\omega)|\xi|^4t}|\hat{\phi}(\xi)|^2d\xi \\
&\qquad
  \leq C\|\phi\|_{L^1_1}^2
  \int_{\mathbb{R}^n}|\xi|^{2(k+1)}e^{-2\gamma(\omega)|\xi|^4t}d\xi
  \leq Ct^{-{n\over4}-{{k+1}\over2}}\|\phi\|_{L^1_1}^2,
\end{split}
\end{equation*}
which proves the first estimate in \eqref{6.15}. 
This completes the proof. 
\end{proof}

\begin{proof}[Proof of Theorem \ref{22}]\  
The function $G_0(x,t+1)$ is a solution of \eqref{limit-eq} with the 
initial data $G_0(x,1)=\phi_0(x)$ in \eqref{2.9}. Therefore it is 
expressed in the form $G_0(x,t+1)=(G_0(t)*\phi_0)(x)$. 
Consequently, we can write 
\begin{equation}\label{6.16}
\begin{split}
u(t)-MG_0(t+1)
&=(u-\bar{u})(t)+\{\bar{u}(t)-G_0(t)*(u_0+u_1)\} \\[1mm]
&+G_0(t)*(u_0+u_1-M\phi_0),
\end{split}
\end{equation}
where $M=\int_{\mathbb{R}^n}(u_0+u_1)(x)dx$. 
The first two terms on the right hand side of \eqref{6.16} have 
already been estimated in Proposition \ref{p51} and Lemma \ref{Su1}, 
respectively. Therefore, for the proof of \eqref{asymp-1} and 
\eqref{asymp-2}, it suffices to show the following decay estimates 
for the last term in \eqref{6.16}. 
\begin{equation}\label{6.17}
\begin{split}
&\|\p^k_xG_0(t)*(u_0+u_1-M\phi_0)\|_{H^{s-k}}
  \leq CE_2(1+t)^{-{n\over8}-{{k+1}\over4}}, \\[1mm]
&\|\p^k_x\p_tG_0(t)*(u_0+u_1-M\phi_0)\|_{H^{s-4-k}}
  \leq CE_2(1+t)^{-{n\over8}-{{k+5}\over4}}.
\end{split}
\end{equation}
To see this, we observe that $\phi_0(x)$ is a rapidly decreasing 
function satisfying $\int_{\mathbb{R}^n}\phi_0(x)dx=1$. 
It then follows that 
$\|\phi_0\|_{H^s}+\|\phi_0\|_{L^1_1}\leq C$ and 
$$
\int_{\mathbb{R}^n}(u_0+u_1-M\phi_0)(x)dx=0.
$$
Therefore, applying \eqref{6.14} and \eqref{6.15} with 
$\phi=u_0+u_1-M\phi_0$, we easily get the desired estimates 
in \eqref{6.17}. Thus the proof of Theorem \ref{22} is complete. 
\end{proof}

\noindent{\bf Acknowledgments.}\ 
This work was partially supported by Grant-in-Aid for JSPS Fellows.



\begin{thebibliography}{}

\bibitem{BL}
M.E. Bradley and S. Lenhart, 
Bilinear spatial control of the velocity term in a Kirchhoff plate 
equation, 
\emph{Electronic J. Differential Equations}, {\bf 2001} (2001), 1-15.

\bibitem{Bu}
C. Buriol, 
Energy decay rates for the Timoshenko system of thermoelastic plates, 
\emph{Nonlinear Analysis}, {\bf 64} (2006), 92-108. 

\bibitem{CBBP}
R.C. Char$\tilde{\rm a}$o, E. Bisognin, V. Bisognin and A.F. Pazoto, 
Asymptotic behavior for a dissipative plate equation in $\R^N$ 
with periodic coefficients, 
\emph{Electronic J. Differential Equations}, {\bf 46} (2008), 23 pp. 

\bibitem{CL}
C.R. da Luz and R.C. Char$\tilde{\rm a}$o, 
Asymptotic properties for a semilinear plate equation in unbounded 
domains, 
\emph{J. Hyperbolic Differential Equations}, 6 (2009), 269-294.

\bibitem{DK}
A.D. Drozdov and V.B. Kolmanovskii, 
Stability in viscoelasticity, 
North-Holland Series in Applied Mathematics and Mechanics, {\bf 38}, 
North-Holland Publishing Co., Amsterdam, 1994.

\bibitem{E}
Y. Enomoto, 
On a thermoelastic plate equation in an exterior domain, 
\emph{Math. Meth. Appl. Sci.}, {\bf 25} (2002), 443-472. 

\bibitem{FL}
M. Fabrizio and B. Lazzari, 
On the existence and the asymptotic stability of solutions for 
linear viscoelastic solids, 
\emph {Arch. Rational Mech. Anal.}, {\bf 116} (1991), 139-152.

\bibitem{HK}
T. Hosono and K. Kawashima, 
Decay property of regularity-loss type and application to some 
nonlinear hyperbolic-elliptic system, 
\emph{Math. Models Meth. Appl. Sci.}, {\bf 16} (2006), 1839-1859.

\bibitem{IK}
K. Ide and S. Kawashima, 
Decay property of regularity-loss type and nonlinear effects for 
dissipative Timoshenko system, 
\emph{Math. Models Meth. Appl. Sci.}, {\bf 18} (2008), 1001-1025.

\bibitem{Le}
H.J. Lee, 
Uniform decay for solution of the plate equation with a boundary 
condition of memory type, 
\emph{Trends in Math.}, {\bf 9} (2006), 51-55.

\bibitem{LW}
Y. Liu and W. Wang, 
The pointwise estimates of solutions for dissipative wave equation 
in multi-dimensions, 
\emph{Discrete and Continuous Dynamical Systems}, 
{\bf 20} (2008), 1013-1028.

\bibitem{LZ1}
Z. Liu and S. Zheng, 
On the exponential stability of linear viscoelasticity and 
thermoviscoelasticity, 
\emph{Quart. Appl. Math.}, {\bf 54} (1996), 21-31.

\bibitem{LZ2}
Z. Liu and S. Zheng, 
Semigroups associated with dissipative systems, 
Chapman $\&$ Hall/CRC, London, 1999.

\bibitem{Ma}
A. Matsumura, 
On the asymptotic behavior of semi-linear wave equations, 
\emph{Publ. Res. Inst. Math. Sci.}, {\bf 12} (1976), 169-189.

\bibitem{MNV}
J.E. Mu$\tilde{\rm n}$oz Rivera, M.G. Naso and F.M. Vegni, 
Asymptotic behavior of the energy for a class of weakly dissipative 
second-order systems with memory, 
\emph{J. Math. Anal. Appl.}, {\bf 286} (2003), 692-704.

\bibitem{Mu}
J.E. Mu$\tilde{\rm n}$oz Rivera, 
Asymptotic behavior in linear viscoelasticity, 
\emph{Quart. Appl. Math.}, {\bf 52} (1994), 628-648.

\bibitem{MZ}
G. Perla Menzala and E. Zuazua, 
Timoshenko's plate equations as a singular limit of the dynamical 
von K\'{a}rm\'{a}n system, 
\emph{J. Math. Pures Appl.}, {\bf 79} (2000), 73-94. 

\bibitem{Ni}
K. Nishihara, 
$L^p-L^q$ estimates of solutions to the damped wave equation in 
$3$-dimensional space and their applications, 
\emph{Math. Z.}, {\bf 244} (2003), 631-649.

\bibitem{Pa}
J.Y. Park, 
Bilinear boundary optimal control of the velocity terms in a 
Kirchhoff plate quation, 
\emph{Trends in Math.}, {\bf 9} (2006), 41-44.

\bibitem{PVM}
A.F. Pazoto, J.C. Vila Bravo and J.E. Mu$\tilde{\rm n}$oz Rivera, 
Asymptotic stability of semigroups associated to linear weak 
dissipative systems, 
\emph{Math. Computer Modeling}, {\bf 40} (2005), 387-392.

\bibitem{RR}
J.E. Mu$\tilde{\rm n}$oz Rivera and R. Racke, 
Global stability for damped Timoshenko systems, 
\emph{Discrete and Continuous Dynamical Systems}, 
{\bf 9} (2003), 1625-1639.

\bibitem{Sa}
J.R. Luyo S\'{a}nchez, 
O sistema din\'{a}mico de von K\'{a}rm\'{a}n en dom\'{i}nios n\'{a}o 
limitados \'{e} globalmente bem posto no sentido de Hadamard: 
An\'{a}lise do seu limite singular, 
Doctoral Thesis, Institute of Mathematics, Federal University 
of Rio de Janeiro, Rio de Janeiro, Brazil, 2003. 

\bibitem{SK}
Y. Sugitani and S. Kawashima, 
Decay estimates of solutions to a semilinear dissipative 
plate equation, 
\emph{J. Hyperbolic Differential Equations}, accepted. 

\bibitem{Te}
R. Teman, 
Navier-Stokes Equations, 
Studies in Mathematics and Its Applications, Vol. 2, Revised Edition, 
North-Holland, Amsterdam, New York, Oxford, 1979. 

\bibitem{TY}
G. Todorova and B. Yordanov, 
Critical exponent for a nonlinear wave equation with damping, 
\emph{J. Differential Equations}, {\bf 174} (2000), 464-489.

\bibitem{ZZ}
X. Zhang and E. Zuazua, 
On the optimality of the observability inequalities for Kirchhoff 
plate systems with potentials in unbounded domains, 
\emph{Hyperbolic Problems: Theory, Numerics and Applications}, 
Springer, 2008, 233-243.

\end{thebibliography}
\end{document}